\documentclass[a4paper,twoside,10pt]{article}
\usepackage[a4paper,left=2cm,right=2cm, top=3cm, bottom=3cm]{geometry}
\usepackage{cuted}

\usepackage{enumerate}
\usepackage{amsthm}
\usepackage{mathrsfs}
\usepackage{mathtools}
\usepackage{accents}
\usepackage{graphicx}
\usepackage{orcidlink}
\usepackage{amscd,amsmath,amssymb,mathrsfs,bbm,listings}
\usepackage{cancel}
\usepackage{tikz}
\usetikzlibrary{decorations.pathreplacing}
\usepackage{epstopdf}
\graphicspath{{Figs/}}
\usepackage{amsmath}
\usepackage{footmisc}
\usepackage{amsthm}
\usepackage{amssymb}
\usepackage{stmaryrd}
\SetSymbolFont{stmry}{bold}{U}{stmry}{m}{n}
\usepackage{float}
\usepackage{bigints}
\usepackage{cite}
\usepackage{color}
\usepackage[abs]{overpic}
\usepackage[font=footnotesize,labelfont=bf]{caption}
\usepackage{cases}
\usepackage{tikz}
\usepackage{algorithmic}
\usepackage{rotating}
\usepackage{blkarray}
\usetikzlibrary{matrix,calc,arrows}
\usepackage{filecontents}
\usepackage{pgfplots}
\usepackage{subcaption}
\usepackage{soul,xcolor}
\usepackage{enumitem}  
\usepackage{verbatim}

\usepackage{graphicx}
\usepackage{alphabeta}
\usepackage{mwe}
\usepackage{algorithm}
\usepackage{array}
\usepackage{multirow}

\newtheorem{theorem}{Theorem}[section]

\newtheorem{remark}[theorem]{Remark}

\newcommand{\eremk}{\hbox{}\hfill\rule{0.8ex}{0.8ex}}

\newcolumntype{C}{>{\centering\arraybackslash}p{2.5cm}}
\numberwithin{equation}{section}
\hypersetup{
    colorlinks=true,
    citecolor=red,
    linkcolor=blue,
    pdfborder={0 0 0}
}

\newcommand{\N}{\mathbb{N}}
\newcommand{\R}{\mathbb{R}}
\newcommand{\bx}{\boldsymbol{x}}
\newcommand{\bn}{\boldsymbol{n}}

\newcommand{\Norm}[2]{\|#1\|_{#2}}

\newcommand{\Ec}{E_c}
\newcommand{\Esigma}{E_{\vsigma}}
\newcommand{\dpt}{\partial_t}

\newcommand{\QT}{Q_T}
\newcommand{\calD}{\mathcal{D}}
\newcommand{\bnOmega}{\bn_{\Omega}}
\newcommand{\bnK}{\bn_{K}}

\newcommand{\D}{\mathbf{D}}

\newcommand{\ba}{\mathbf{a}}
\newcommand{\daxn}{\text{d}_{\mathrm{axn}}}
\newcommand{\dext}{\text{d}_{\mathrm{ext}}}

\newcommand{\vz}{\boldsymbol{z}}

\newcommand{\vsigma}{\boldsymbol{\sigma}}

\newcommand{\vr}{\boldsymbol{r}}

\newcommand{\In}{I_n}
\newcommand{\tn}{t_n}
\newcommand{\tnmo}{t_{n - 1}}

\newcommand{\Th}{\mathcal{T}_h}
\newcommand{\Tt}{\mathcal{T}_{\tau}}

\newcommand{\npo}{(n + 1)}

\newcommand{\LDG}{{\mathrm{LDG}}}
\newcommand{\nablaLDG}{\nabla_{\LDG}}
\newcommand{\divLDG}{\mathrm{div}_{\LDG}}
\newcommand{\DG}{{\mathrm{DG}}}

\newcommand{\dx}{\, \mathrm{d}\boldsymbol{x}}

\newcommand{\dt}{\, \mathrm{d}t}

\newcommand{\Wh}{\mathcal{W}^{\,\ell}(\Th)}
\newcommand{\Rh}{\boldsymbol{\mathcal{R}}^\ell(\Th)}
\newcommand{\vh}{v_h}
\newcommand{\wh}{w_{h}}

\newcommand{\whp}{w_{p, h}}
\newcommand{\whq}{w_{q, h}}
\newcommand{\zh}{\boldsymbol{z}_{h}}

\newcommand{\sigmah}{\boldsymbol{\sigma}_{h}}

\newcommand{\rh}{\boldsymbol{r}_{h}}

\newcommand{\phih}{\boldsymbol{\phi}_h}

\newcommand{\psih}{\psi_h}

\newcommand{\cvW}{\mathbf{W}_{h}}

\newcommand{\cvZ}{\mathbf{Z}_{h}}

\newcommand{\cvSigma}{\mathbf{\Sigma}_{h}}

\newcommand{\cvR}{\mathbf{R}_{h}}

\newcommand{\MM}{\mathbb{M}}

\newcommand{\GG}{\mathcal{G}}

\newcommand{\diffsigma}{\mathcal{D}_{\mathbf{\Sigma}}\!}
\newcommand{\diffw}{\mathcal{D}_{\mathbf{W}}\!}

\newcommand{\PiW}{\Pi_{\mathcal{W}}}
\newcommand{\PiR}{\boldsymbol{\Pi}_{\boldsymbol{\mathcal{R}}}}
\newcommand{\Pp}[2]{\mathbb{P}^{#1}(#2)}

\def\Fh{\mathcal{F}_h}
\def\Fho{\mathcal{F}_h^\mathcal{I}}
\def\FhN{\mathcal{F}_h^{\mathcal{N}}}

\newcommand{\mvl}[1]{\{ \!\!\{#1\}\!\!\}}  
\newcommand{\jump}[1]{\llbracket #1\rrbracket}

\newcommand{\mvlparamF}{\gamma_F}

\def\bdfp{\nu\,\!}
\newcommand{\BDF}[1]{\textrm{BDF}{#1}}

\title{A structure-preserving LDG discretization of the Fisher-Kolmogorov equation for modeling neurodegenerative diseases}

\author{Paola F. Antonietti\,\orcidlink{0000-0002-2138-3878}\thanks{MOX-Dipartimento di Matematica, Politecnico di Milano, Piazza Leonardo da Vinci 32, Milan, 20133, Italy (\href{mailto:paola.antonietti@polimi.it}{paola.antonietti@polimi.it}, \href{mailto:mattia.corti@polimi.it}{mattia.corti@polimi.it})} 
\and
Mattia Corti\,\orcidlink{0000-0002-7014-972X}\footnotemark[1]
\and
Sergio G\'omez\,\orcidlink{0000-0001-9156-5135}\thanks{Department of Mathematics and Applications, University of Milano-Bicocca, Via Cozzi 55, Milan, 20125, Italy (\href{mailto:sergio.gomezmacias@unimib.it}{sergio.gomezmacias@unimib.it})}\ \thanks{IMATI CNR ``E. Magenes", Via Ferrata 5, 27100, Pavia, Italy} 
\and 
Ilaria Perugia\,\orcidlink{0000-0003-1368-2883}\thanks{Faculty of Mathematics, University of Vienna, Oskar-Morgenstern-Platz 1, 1090 Vienna, Austria (\href{mailto:ilaria.perugia@univie.ac.at}{ilaria.perugia@univie.ac.at})}
}

\date{}

\begin{document}
\maketitle

\begin{abstract}
\noindent 
This work presents a structure-preserving, high-order, unconditionally stable numerical
method for approximating the solution to the Fisher-Kolmogorov equation on polytopic meshes, with a particular focus on its application in simulating misfolded protein spreading in neurodegenerative diseases. 
The model problem is reformulated using an entropy variable to guarantee solution positivity, boundedness, and satisfaction of a discrete entropy-stability inequality at the numerical level. The scheme combines a local discontinuous Galerkin 
method on polytopal meshes for the space discretization with a~$ν$-step backward differentiation formula 
for the time integration. Implementation details are discussed, including a detailed derivation of the linear systems arising from Newton’s iteration. 
The accuracy and robustness of the proposed 
method are demonstrated through extensive numerical tests. Finally, the method’s practical performance is demonstrated through simulations of \textalpha{}-synuclein propagation in a two-dimensional brain geometry segmented from MRI data, providing a relevant computational framework for modeling synucleopathies (such as Parkinson’s disease) and, more generally, neurodegenerative diseases. 
\end{abstract}


\section{Introduction}

The Fisher-Kolmogorov (FK) equation, also known as the Fisher-Kolmogorov-Petrovsky-Piskunov (Fisher-KPP) equation, is a reaction-diffusion equation that models the dynamics of biological populations or advantageous genes \cite{1937_Fisher}. The FK equation combines a diffusion term with a nonlinear reaction component, accounting for growth and saturation effects. In recent years, its application has become increasingly prominent in modeling neurodegenerative diseases, particularly within the framework of the \emph{prion-like} hypothesis \cite{fornari_prion-like_2019}. This hypothesis suggests that misfolded proteins, such as \textalpha{}-synuclein in Parkinson's disease or dementia with Lewy Bodies (DLB) \cite{beach_unified_2009}, and amyloid-\textbeta{} and tau in the context of Alzheimer's disease \cite{bloom2014amyloid}, propagate throughout the brain like infectious agents, transitioning from one region to another.
Consequently, the FK equation provides a valuable mathematical framework capable of addressing misfolded proteins dynamics, local proliferation dynamics, and their spread through brain networks. 
At the continuous level, the key variables—population density or pathological protein concentration—are inherently nonnegative and constrained by physical and biological limits. Indeed, from a mathematical viewpoint, it can be proved that the solution to the FK equation with homogeneous Neumann boundary conditions remains nonnegative and bounded as long as the initial conditions are nonnegative and the reaction term satisfies appropriate biological constraints corresponding to the absence and saturation of the population or pathological agent, respectively \cite{salsa_partial_2022}.  When moving from the infinite-dimensional setting to the discrete one, it is crucial to 
develop numerical discretization methods that 
preserve the nonnegativity and boundedness of the discrete solution. Additionally, they should reproduce key structural properties of the continuous problem to ensure consistency with the underlying model. Such schemes are 
referred to as \emph{structure-preserving} (encompassing positivity- and bound-preservation). 

Discontinuous Galerkin (DG) methods have been effectively employed to discretize FK equations owing to their flexibility in handling complex and heterogeneous computational domains, as well as their ability to accommodate (element-wise) high-order approximations, which is pivotal for accurately capturing traveling-wave solutions. These features are particularly significant in the modeling of neurodegenerative diseases, where the underlying geometry is complex, with intricate details and highly heterogeneous media, and where phenomena occur at long time scales.
DG methods are now well established; for a detailed account of their historical development, we refer the reader to~\cite{2002_Arnold_Brezzi_Cockburn_Marini}.
To further enhance geometric flexibility, DG methods on polytopal grids (PolyDG) have been proposed for the numerical discretization of a broad range of differential problems, see e.g., \cite{2009_Antonietti_Brezzi_Marini,Bassi2012,Antonietti2013,Cangiani2014}, the references therein, 
the comprehensive monograph~\cite{Cangiani_Dong_Georgoulis:2017}, 
and the review~\cite{Antonietti2021Review}.
PolyDG discretizations of FK equations have been studied in~\cite{Corti_Bonizzoni_Dede_Quarteroni_Antonietti:2023} and have been extended to encompass more physics-informed models in neurodegeneration, including the heterodimer model~\cite{Antonietti_Bonizzoni_Corti_Dallolio:2024,2024_Antonietti_Corti_Lorenzon}. Nonetheless, numerical discretizations written in the primal variable, i.e.~the pathological protein concentration for our application of interest, are not inherently structure-preserving. 
In some cases, 
they may produce discrete solutions with nonphysical negative values or values that exceed 
biologically plausible limits. 
This can occur, for example, when approximating steep gradients or using coarse meshes. 

The design of structure-preserving DG methods is a 
well-established area of research.
Positivity-preserving DG schemes for parabolic equations have been developed, e.g., in~\cite{2015_Guo_Yang,2018_Sun_Carrillo_Shu,2016_Zhang_et_al}. 
For the FK equation, positivity-preserving schemes have been explored in the context of finite-difference approximations (see, e.g.,  \cite{2012_Macias_et_al,2017_Hasnain_Saqiband_Suleiman_Mashat}) and of continuous finite-element discretizations (see, e.g., \cite{2017_Yadav_Jiwari}).
A possible approach to devise a positivity-preserving scheme for the FK equation involves rewriting the original model by changing variables. Instead of solving for the original field, the system is reformulated 
in terms of an 
entropy variable. 
In this context, an entropy-stable, structure-preserving DG scheme on simplicial grids for the FK equation has been proposed and analyzed in~\cite{Bonizzoni_Braukhoff_Jungel_Perugia:2020} and further extended to polytopal grids in~\cite{Corti_Bonizzoni_Dede_Quarteroni_Antonietti:2023}. 
The latter work also provided a thorough
validation in the context of modeling the accumulation and spread of prionic proteins, showing accurate solutions for biologically meaningful, stable, and accurate simulations of disease progression. In~\cite{Gomez_Jungel_Perugia:2024}, a structure-preserving local discontinuous Galerkin (LDG) method on simplicial grids for nonlinear cross-diffusion systems 
that exploits the boundedness-by-entropy framework of~\cite{Jungel:2015} is presented and analyzed. 

The LDG method was first introduced in~\cite{CockburnShu1998} for time-dependent convection-diffusion systems. For problems
involving second-order elliptic operators, the classical interior penalty DG (IPDG) method~\cite{1982_Arnold} can be derived directly
by integrating by parts, symmetrizing, and adding suitable stability terms. In contrast, the LDG method reformulates the PDE by introducing an auxiliary variable, which allows the original second-order equation to be written as a system of first-order equations. The numerical fluxes for the auxiliary vector-valued variable are chosen as (weighted) averages of the traces of the auxiliary variable itself (instead of the traces of the gradient of the primal variable, as in IPDG).
As a result, the interface terms do not depend explicitly on differential operators. The auxiliary variable can then be eliminated locally at the element level, avoiding an increase in the number
of unknowns. For problems with nonlinear diffusion, as those considered in ~\cite{Gomez_Jungel_Perugia:2024,Zheng_Xu:2024}, the strategy of introducing auxiliary unknowns can be further extended to reformulate the problem so that nonlinearities do not appear under differential operators. This ensures that the penalty parameters remain free of any nonlinear dependence and that the equations containing nonlinearities lead to block-diagonal structures in the system matrix.

In the context of linear diffusion, the variable transformation utilized in \cite{Bonizzoni_Braukhoff_Jungel_Perugia:2020}, see also \cite{Corti_Bonizzoni_Dede_Quarteroni_Antonietti:2023}, to enforce positivity of solutions introduces a nonlinearity in the diffusion tensor, which complicates the practical implementation of the IPDG discretization proposed there. First, the penalty function in the IPDG discretization of the diffusion term depends on the new auxiliary variable obtained from the exponential transformation and necessitates careful calibration. Furthermore, due to the fact that 
nonlinearities appear in the differential operator, Newton’s iterations may exhibit poor conditioning, especially when approximation spaces with large polynomial degrees are employed.

This paper aims to specialize 
the approach proposed in \cite{Gomez_Jungel_Perugia:2024} 
to the FK equation. The main idea is to exploit the LDG paradigm, as well as the underlying entropy structure of the problem, and first rewrite it 
in appropriately chosen auxiliary variables in such a way
that the nonlinearities do not occur within differential operators or interface terms, but only in 
zero-order terms. At the same time, such a transformation allows for the physical positivity and boundedness 
to be satisfied in a strong sense. Next, the 
LDG method on polytopal grids for the space 
discretization, combined with a~$\nu$-step backward differentiation formula (BDF) for the time discretization, 
is employed for the numerical approximation. 
The resulting method inherently preserves positivity and boundedness of solutions for any time discretization, as those properties are guaranteed by the reformulation in the new variables. Furthermore, in the case of backward Euler (BDF1), we prove that the fully discrete problem satisfies an entropy stability inequality similar to that of the original problem. We chose the BDF family of methods for high-order time integration because, compared to the backward Euler method, they require modifications only to the right-hand side, without affecting the computation of the matrices in the tangent problem during Newton's iterations. Nonetheless, we point out that the stability analysis of the fully discrete problem with high-order time integration remains an open issue. Additionally,
since nonlinearities arise only in volume terms, these are treated independently cell-wise, exploiting the fact that volume terms are block diagonal in the DG framework. 
We also emphasize that, since the nonlinearities do not appear within interface terms, the DG penalty parameters can be chosen exactly as in the linear case, without any nonlinear dependence in their definition. Moreover, the Newton solver benefits from the block-diagonal structure corresponding to the equations containing the nonlinearities, which can be treated independently cell-wise.

The remaining part of the manuscript is structured as follows. In Section~\ref{sec:model_problem_numerical_method}, we introduce the model problem, reformulate it using a suitable entropy variable, and present its structure-preserving LDG space discretization on polytopal grids, coupled with a~$\nu$-step BDF time integration scheme. Section~\ref{sec:model_problem_numerical_method} also discusses implementation details, including an explicit derivation of the linear systems stemming from Newton’s iteration. Section~\ref{sec:verification} provides 
an extensive numerical verification of the accuracy of the proposed structure-preserving BDF-LDG method through test cases with known analytical solutions. Finally, Section~\ref{sec:brain} presents computations aimed at demonstrating the practical performance of the proposed method in modeling \textalpha{}-synuclein protein spreading, which is involved in relevant neurodegenerative diseases such as Parkinson’s disease, dementia with Lewy bodies, and Alzheimer’s disease with Lewy bodies. The results in this section have been obtained on a two-dimensional section of a brain geometry, segmented from a structural
magnetic resonance image.

\section{Model problem and numerical discretization }\label{sec:model_problem_numerical_method}

Let~$\Omega \subset \R^d$ ($d\in \{2, 3\}$) be a polygonal/polyhedral space domain with Lipschitz boundary~$\Gamma:=\partial \Omega$ and outward-pointing normal unit vector~$\bnOmega$, and let~$(0,T)$ be a time interval, with~$T>0$. 
In the space--time cylinder~$\QT := \Omega \times (0, T)$, we consider the following initial and boundary value problem for the Fisher-Kolmogorov (FK) equation: find~$c:\QT \rightarrow \R$ such that 
\begin{subequations}
\label{EQN::FISHER-KPP}
\begin{alignat}{3}
\dpt c - \nabla \cdot (\D \nabla c) & = \alpha\, c (1 - c)=:f(c) & & \quad \text{ in } \QT,\\
(\D \nabla c) \cdot \bnOmega & = \mathbf{0}  & & \quad  \text{ on } \Gamma \times (0, T),\\
c(\cdot, 0) & = c_0 & & \quad \text{ in } \Omega.
\end{alignat}
\end{subequations}
Here, the coefficient~$\alpha = \alpha(\bx) \in L^{\infty}(\Omega)$ is strictly positive, and the diffusion tensor~$\D=\D(\bx)\in L^\infty(\Omega)^{d\times d}$ is assumed to be symmetric and uniformly positive definite, namely there exists a constant~$D_0>0$ such that
\begin{equation}
\label{EQ:POSITIVE-DEFINITE-DIFFUSION}
\vz^\top \D\vz
\ge D_0\,|\vz|^2 \qquad \forall\vz\in\R^d.
\end{equation}
Additionally, we assume that~$0\le c_0\le 1$ a.e.~in~$\Omega$, so that
\begin{equation}\label{EQ::cbound}
0\le c(\bx,t) \le 1\qquad \text{for almost all~$(\bx,t)\in Q_T$};
\end{equation}
see~\cite[\S2.10.2]{salsa_partial_2022}.

Next, we introduce our structure-preserving LDG method. We start by reformulating the problem in an entropy variable in Section~\ref{SEC::REFORMULATION}, then we define the method in Section~\ref{SEC::METHOD}, and write it in matrix form in Section~\ref{SEC::MATRIX}.
Henceforth, given a domain~$\calD$, we use~$(\cdot, \cdot)_{\calD}$ to denote the~$L^2(\calD)$ inner product, and~$|\calD|$ its measure.

\subsection{Reformulation in terms of an
entropy variable}\label{SEC::REFORMULATION}

The change of variable
\begin{equation}\label{EQ::UPUQ}
c=u(w):=\frac{e^{w}}{1+e^{w}}, \quad w: \QT\to \R,
\end{equation}
enforces that~$c$ satisfies the bound in~\eqref{EQ::cbound}. In the setting of~\cite{Jungel:2015,Gomez_Jungel_Perugia:2024}, this corresponds to writing~$u=(s')^{-1}$, with the \emph{entropy function}~$s$ defined by
\begin{equation*}
s(c):=c\log c+\left(1-c\right)\log\left(1-c\right) + \log 2 \ge 0, 
\end{equation*}
for which
\[
s'(c)=\log c -\log\left(1-c\right)=\log\left(\frac{c}{1-c}\right),
\qquad
s''(c)=\frac{1}{c(1-c)}.
\] 
Such an entropy function satisfies the following properties:
\begin{subequations}
\label{EQN::PROPERTIES}
\begin{alignat}{3}
\label{EQN::PROPERTIES-A}
s''(c) & \geq 4
& & \quad \forall c \in (0, 1), \\
\label{EQN::PROPERTIES-B}
|\alpha(\bx)  c(1 - c) s'(c)| & \le C_f \Norm{\alpha}{L^{\infty}(\Omega)}  & & \quad \text{for a.e.}\ \bx \in \Omega, \ \forall c \in (0, 1),
\end{alignat}
\end{subequations}
with~$C_f = \max_{c \in [0, 1]} |c(1 - c) s'(c)|\le 0.25$. 
As~$w=s'(c)$, the chain rule gives
\begin{equation}\label{EQ::CHAINRULE}
\nabla w=s''(c)\nabla c.
\end{equation}

We introduce the following auxiliary variables in~$Q_T$:
\begin{subequations}
\label{EQ::HETERODIMER_VARIABLES}
\begin{align}
\label{EQ::var_c}
w\ \, \text{s.t.}\ \, c &\,= u(w), \\
\label{EQ::Z}
\vz  &:= -\nabla w, \\
\label{EQ::SIGMA}
    \D s''(c) \, \vsigma  &:= - \D s''(c)  \nabla c = \D \vz,  \\
    \label{EQ::Q}
    \vr  &:= \D \vsigma.
\end{align}
\end{subequations}
The first identity in equation~\eqref{EQ::SIGMA} imposes~$\vsigma=-\nabla c$, due to the invertibility of $\D$ and the fact that~$s''$ is bounded away from zero. More precisely,~$s''(c)\D$ is uniformly positive definite with constant~$4 D_0$, where $D_0$ is the constant in~\eqref{EQ:POSITIVE-DEFINITE-DIFFUSION}:
\begin{equation}\label{EQ::UPD}
\vz^\top \left(s''(c)\D\right)\vz
\ge\left(\inf_{(\bx,t)\in Q_T}s''\left(c(\bx,t)\right)\right)
\,D_0\,|\vz|^2\ge 4 D_0\,|\vz|^2 \qquad \forall\vz\in\R^d.
\end{equation}
After applying the change of variables in equation~\eqref{EQ::var_c}, the definition of the auxiliary variables follows the standard LDG approach, avoiding the presence of non-linearities under differential operators. This is achieved by imposing the chain rule in equation~\eqref{EQ::CHAINRULE}. Specifically, we define~$\vr:=-\D\nabla c$ (equation~\eqref{EQ::Q}), writing~$\vsigma=-\nabla c$ (first identity in~\eqref{EQ::SIGMA}). 
The second identity in~\eqref{EQ::SIGMA}, along with equation~\eqref{EQ::Z}, enforces the chain rule in equation~\eqref{EQ::CHAINRULE}, avoiding that~$\nabla c=\nabla\left(u(w)\right)$ appears in the formulation.
The variable~$\vr$ then represent the flux of~$c$ and will appear explicitly in the formulation (see~\eqref{EQN::REWRITTEN} below). 

Problem~\eqref{EQN::FISHER-KPP} is reformulated in the variables defined in~\eqref{EQ::HETERODIMER_VARIABLES} as follows: find~$w: \QT \rightarrow \R$ and~$\vr: \QT \rightarrow \R^d$ such that
\begin{subequations}
\label{EQN::REWRITTEN}
\begin{alignat}{3}
\dpt c + \nabla \cdot \vr & = \alpha\, c (1 - c)=f(c) & & \quad \text{ in } \QT,\\
\vr \cdot \bnOmega & = 0  & & \quad  \text{ on } \Gamma \times (0, T), \\
c(\cdot, 0) & = c_{0}  & & \quad \text{ in } \Omega\times\{0\},
\end{alignat}
\end{subequations}
where, in~$Q_T=\Omega\times (0,T)$,~$c=u(w)$ is understood. 
The change of variable~\eqref{EQ::UPUQ} and the reformulation~\eqref{EQN::REWRITTEN} of the FK problem with variables defined in~\eqref{EQ::HETERODIMER_VARIABLES} are at the basis of the structure-preserving discretization we introduce in the following section.

We conclude this section by noting that the solution~$c$ of the FK problem satisfies the following \emph{entropy stability estimate}, which follows from testing the variational formulation of~\eqref{EQN::FISHER-KPP} with~$w$, and using~\eqref{EQN::PROPERTIES-B}, \eqref{EQ::UPD}, and the chain rule in~\eqref{EQ::CHAINRULE}:
\begin{equation}
\label{EQN::CONTINUOUS-ENTROPY-STABILITY}
    \int_{\Omega} s(c(\cdot, \tau)) \dx + 4 D_0 \int_0^{\tau} \Norm{\nabla c}{L^2(\Omega)^d}^2 \dt \leq \int_{\Omega} s(c_0) \dx + C_f \Norm{\alpha}{L^{\infty}(\Omega)}\tau |\Omega| \qquad \text{ for all } 0 < \tau \leq T.
\end{equation}

\subsection{The structure-preserving BDF-LDG method}\label{SEC::METHOD}
In this section, we present a numerical discretization of the FK problem~\eqref{EQN::FISHER-KPP} based on formulation~\eqref{EQN::REWRITTEN} with auxiliary variables defined in~\eqref{EQ::HETERODIMER_VARIABLES}.

\bigskip
\paragraph{Meshes.}
We partition the space domain~$\Omega$ by a locally quasi-uniform polytopal mesh~$\Th$ with mesh size~$h:=\max_{K\in\Th}h_K$, where~$h_K$ denotes the diameter of the element~$K$.
We also partition the time interval~$(0, T)$ by a mesh~$\Tt$ defined by the points~$0 := t_0 < t_1 < \ldots < t_N := T$. For~$n = 1, \ldots, N$, we define the time interval~$\In := (\tnmo, \tn)$ and the time step~$\tau_n := \tn - \tnmo$. 

As usual in the DG setting, we denote the set of all the mesh facets in~$\Th$ by~$\Fh = \Fho \cup \FhN$, where~$\Fho$ and~$\FhN$ are the sets of internal and (Neumann) boundary facets, respectively. Additionally, we define the mesh-size function~$\mathsf{h} \in L^{\infty}(\Fho)$ as
\begin{equation}
\label{EQN::DEF-h}
\mathsf{h}(\bx) := \eta_F^{-1} \left[\dfrac{1}{2}\left(\left(\dfrac{|K_1|}{m_{K_1}|F|}\right)^{\theta} + \left(\dfrac{|K_2|}{m_{K_2}|F|}\right)^{\theta}\right)\right]^{1/\theta} \quad \text{ if }\bx \in F, \text{ and~$F\in \Fho$ is shared by~$K_1, K_2\in \Th$},
\end{equation}
where~$\eta_F>0$ is a parameter independent of the mesh size, which is chosen in~\eqref{eq:param} below in terms of the diffusion tensor and the polynomial degree in the space discretization, $\theta$ is the power-mean exponent, and $m_{K_*}$ is the number of facets of the polytopal element~$K_*$. 
In the numerical experiments presented in Section~\ref{sec:verification}, we use~$\theta=-1$, which corresponds to the harmonic mean, since this classical choice works well on families of regular Voronoi meshes, as those used there. On the other hand, inspired by~\cite[\S3 and~\S4]{Dong_Georgoulis:2022}, where more general polytopal meshes are considered, we set $ \theta=1/2$ in the experiments of Section~\ref{sec:brain}, which, due to the complicated geometry, are carried out on an agglomerated mesh.

For each element~$K\in\Th$, let~$\bnK$ denote the unit normal vector in~$d$ dimensions to~$\partial K$, pointing away from~$K$. Given any piecewise smooth, scalar function~$\mu$ and any~$d$-vector-valued function~$\boldsymbol{\mu}$, we define the normal jumps and weighted mean values as follows: 
on each facet $F\in \Fho$ shared by two elements~$K_1$ and~$K_2$ of the mesh~$\Th$, we set
\begin{alignat*}{3}
\jump{{\mu}}_{\sf N} & := {{\mu}}_{|_{K_1}} \bn_{K_1} + 
{{\mu}}_{|_{K_2}} \bn_{K_2}, & \quad 
\mvl{\mu}_{\mvlparamF}  & :=  (1 - \mvlparamF) \mu_{|_{K_1}} + \mvlparamF \mu_{|_{K_2}},\\
\jump{\boldsymbol{\mu}}_{\sf N} & := {\boldsymbol{\mu}}_{|_{K_1}} \cdot \bn_{K_1} + 
{\boldsymbol{\mu}}_{|_{K_2}} \cdot \bn_{K_2},
& \quad
\mvl{\boldsymbol{\mu}}_{\mvlparamF} & :=  (1 - \mvlparamF) \boldsymbol{\mu}_{|_{K_1}} + \mvlparamF \boldsymbol{\mu}_{|_{K_2}},
\end{alignat*}
with~$\mvlparamF\in [0,1]$. 

For the definition of~$\eta_F$ and~$\gamma_F$, we follow~\cite[\S2]{Ern2009} 
and set 
\begin{equation}\label{eq:param}
\eta_F = \eta_0 \ell^2  \dfrac{2(\boldsymbol{n}_{K_1}^T\mathbf{D}_{|_{K_1}}\boldsymbol{n}_{K_1})(\boldsymbol{n}_{K_2}^T\mathbf{D}_{|_{K_2}}\boldsymbol{n}_{K_2})}{\boldsymbol{n}_{K_1}^T\mathbf{D}_{|_{K_1}}\boldsymbol{n}_{K_1}+\boldsymbol{n}_{K_2}^T\mathbf{D}_{|_{K_2}}\boldsymbol{n}_{K_2}},\qquad
\mvlparamF=\dfrac{\boldsymbol{n}_{K_1}^T\mathbf{D}_{|_{K_1}}\boldsymbol{n}_{K_1}}{\boldsymbol{n}_{K_1}^T\mathbf{D}_{|_{K_1}}\boldsymbol{n}_{K_1}+\boldsymbol{n}_{K_2}^T\mathbf{D}_{|_{K_2}}\boldsymbol{n}_{K_2}},
\end{equation}
where~$\ell$ is the polynomial degree in the space discretization, and~$\eta_0>0$ is a constant independent of the problem coefficients and discretization parameters. The definitions of the parameters~$\eta_F$ and~$\gamma_F$ are guided by the literature on robust discontinuous Galerkin methods for diffusion problems with strongly heterogeneous coefficients. In these cases, the use of suitable weighted averages and penalty parameters makes the method robust, as shown, e.g., in \cite{Dong_Georgoulis:2022,Ern2009}.
\bigskip

\paragraph{DG spaces for the discretization in space.}
Given a polynomial degree~$\ell \in \N$ with~$\ell \geq 1$, we define the following spaces:
\begin{equation*}
\Wh := \prod_{K \in \Th} \Pp{\ell}{K} \qquad \text{ and } \qquad \Rh := \prod_{K \in \Th} \Pp{\ell}{K}^d,
\end{equation*}
where~$\Pp{\ell}{K}$ denotes the space of scalar polynomials of degree at most~$\ell$ defined on~$K$. Additionally, we denote by~$\PiW$ and~$\PiR$ the~$L^2(\Omega)$- and~$L^2(\Omega)^d$-orthogonal projections in~$\Wh$ and~$\Rh$, respectively.

\bigskip
\noindent
{\bf Discrete differential operators in space.}
We introduce the discrete space gradient and divergence operators. The LDG gradient operator~$\nablaLDG : \Wh \to \Rh$ is defined by 
\begin{equation*}
\left(\nablaLDG \vh , \phih \right)_\Omega
=
\left(\nabla_h \vh - \mathcal{L}(\vh),  \phih \right)_\Omega
\quad\forall  \phih \in \Rh,
\end{equation*}
where~$\nabla_h$ denotes the piecewise gradient operator, and the jump lifting operator $\mathcal{L}:
\Wh \to \Rh$ is given by
\[
\left(\mathcal{L}(\vh), \phih \right)_\Omega
=\sum_{F\in\Fho} \left(\jump{\vh}_{\sf N}, \mvl{\phih}_{1-\mvlparamF}\right)_F
\quad\forall  \phih \in \Rh.
\]
The LDG divergence operator~$\divLDG: \Rh \to \Wh$ is defined by
\begin{equation*}
\left(\divLDG \rh, \psih\right)_\Omega
=
- \left(\rh, \nablaLDG \psih\right)_{\Omega}
\quad\forall  \psih \in \Wh.
\end{equation*}
The definitions of the LDG gradient and divergence operators correspond to choosing the numerical fluxes for~$\wh$ and~$\rh$ in terms of~$\mvl{\wh}_{\gamma_F}$ and~$\mvl{\rh}_{1-\gamma_F}$, respectively. This choice ensures the symmetry of the bilinear form discretizing the diffusion operator, and can be found, for instance, in~\cite[\S3.2]{Castillo_Sequeira:2013}.

\paragraph{BDF time stepping.}
We discretize in time with the BDF with~$\bdfp$ steps. As for~$\bdfp\ge 7$ the BDF is unstable~(see, e.g.,
\cite[Thm.~3.4 in~{\S}III.3]{Hairer_Norsett_Wanner:1993}), we limit ourselves to~$\bdfp\le 6$. The notation is as follows: for a nonautonomous ordinary differential equation~$y'=g\left(y,t\right)$, the $\bdfp$-step BDF scheme, denoted by BDF$\nu$, reads
\[
y^{(n+1)}-\sum_{j=1}^\bdfp a_j(\bdfp)\, y^{(n+1-j)}=\tau_{n+1}\beta(\bdfp)\, g\big(y^{(n+1)},t_{n+1}\big).
\]
From here on, we omit the dependence on~$\bdfp$ and simply denote the coefficients~$\beta$ and~$a_j$. For uniform meshes, the coefficients for~$\nu$ up to 6 are provided, e.g., in~\cite[\S3.12]{Lambert:1991}.
In particular, when~$\bdfp=1$, the BDF method coincides with the backward Euler method, regardless of the mesh.
\paragraph{The BDF-LDG method.}
Given a penalty parameter~$\varepsilon > 0$, our BDF-LDG discretization of the reformulated FK problem~\eqref{EQN::REWRITTEN}, after an initialization phase (e.g., using a one-step method, \cite[{\S}III.1]{Hairer_Norsett_Wanner:1993}), 
is as follows: for~$n = \bdfp-1, \ldots, N - 1$, find~$\wh^{\npo}\in \Wh$ and~$\zh^{\npo}$, $\sigmah^{\npo}$, $\rh^{\npo}\in \Rh$ 
(for brevity, we omit the dependence on~$\varepsilon$), such that
\begin{subequations}
\label{EQ::VARIATIONAL}
\begin{alignat}{3}
\zh^{\npo} &= -\nablaLDG \wh^{\npo}, \\
\label{EQ::VARIATIONAL_SIGMA}
\left(\D s''\big(u(\wh^{\npo})\big)  \sigmah^{\npo},\, \phih\right)_{\Omega} & = \left(\D \zh^{\npo},\, \phih\right)_{\Omega} & & 
\qquad \forall \phih \in \Rh, \\
\rh^{\npo} &= \PiR \big(\D \sigmah^{\npo}\big),\\
\nonumber
\varepsilon \big(\wh^{\npo},\, \psih \big)_{\LDG} + \frac{1}{\tau_{n + 1}\beta} \Big(u(\wh^{\npo}) & - \sum_{j=1}^\bdfp a_ju_h^{(n+1-j)},\, \psih \Big)_{\Omega}  &  \\
\nonumber
+ \left(\divLDG \rh^{\npo},\, \psih\right)_{\Omega} 
+  \sum_{F \in \Fho} \big(& \mathsf{h}^{-1} \jump{\wh^{(n + 1)}}_{\sf N}, \jump{\psih}_{\sf N} \big)_{F} \\
\label{EQ::VARIATIONAL_EQ}
& = \left(f\big(u(\wh^{\npo})\big),\, \psih \right)_{\Omega} & & 
\qquad \forall \psih \in \Wh.
\end{alignat}
\end{subequations}
In~\eqref{EQ::VARIATIONAL_EQ}, $(\cdot, \cdot)_{\LDG}$ is the bilinear form
\[
(w,\psi)_{\LDG} := (\alpha\,w, \psi)_\Omega
+(\D\nablaLDG w,\nablaLDG\psi)_\Omega
+\sum_{F \in \Fho} \big( \mathsf{h}^{-1} \jump{w}_{\sf N}, \jump{\psi}_{\sf N} \big)_{F}.
\]
which is coercive in the norm
\begin{equation*}
\Norm{w}{\DG}^2:=
\Norm{\alpha^{\frac12} w}{L^2(\Omega)}^2
+\Norm{\D^{\frac12}\nabla_h w}{L^2(\Omega)^d}^2
+ \sum_{F \in \Fho}\Norm{{\mathsf{h}}^{-\frac12} \jump{w}_{\sf N}}{L^2(F)^{d}}^2.
\end{equation*}
Furthermore,
the terms~$u_h^{(\kappa)}$ under the summation  symbol, which are transmitted from the previous time steps, are defined as follows:
\begin{alignat*}{6}
u_h^{(\kappa)} & := \begin{cases}
\PiW c_0 & \text{ if } \kappa = 0, \\
\PiW u(\wh^{(\kappa)}) & \text{ if } \kappa > 0,
\end{cases}
\end{alignat*}
namely, $u_h^{(0)}$ is 
the~$L^2(\Omega)$ projection into~$\Wh$ of the initial condition for the \emph{original} variable~$c$; at the subsequent time steps, they are the~$L^2(\Omega)$ projections of the \emph{transformed} variables computed at previous time steps.

As $\D s''\left(u(w)\right)$ is uniformly positive definite, see~\eqref{EQ::UPD}, given~$\wh^{\npo}$ and~$\zh^{\npo}$, equation~\eqref{EQ::VARIATIONAL_SIGMA} determines~$\sigmah^{\npo}$ in a unique way.
\begin{remark}[Role of the penalty term]
\label{RMK:PENALTY-TERM}
The penalty term with parameter~$\varepsilon > 0$ in~\eqref{EQ::VARIATIONAL_EQ} is introduced to prevent~$u(w)$ from getting too close to the extreme values~$0$ and~$1$, as the term~$s''(\cdot)$ in~\eqref{EQ::VARIATIONAL_SIGMA} is singular at those values. Consequently, the penalty term improves the stability and convergence of non-linear solvers such as Newton's method described in Section~\ref{SEC::NEWTON}.
\eremk
\end{remark}

\begin{remark}[Discrete entropy stability]\label{RMK:DISCR_STAB}
For the lowest-order time stepping BDF1, which coincides with the backward Euler scheme, the fully discrete method~\eqref{EQ::VARIATIONAL} satisfies the following discrete version of the continuous entropy stability estimate~\eqref{EQN::CONTINUOUS-ENTROPY-STABILITY}, see~\cite[Thm.~3.1]{Gomez_Jungel_Perugia:2024}:
\begin{equation*}
\begin{split}
\varepsilon \sum_{n = 0}^{N - 1} \tau_{n + 1} \Norm{w_h^{(n + 1)}}{\DG}^2 & + \int_{\Omega} s(u(\wh^{(N)})) \dx  + 4 D_0 \sum_{n = 0}^{N - 1} \tau_{n + 1}\Norm{\sigmah^{(n + 1)}}{L^2(\Omega)^d}^2 \\
& + \sum_{n = 0}^{N - 1} \sum_{F \in \Fho} \tau_{n + 1} \Norm{\mathsf{h}^{-\frac12} \jump{\wh^{(n + 1)}}_{\sf{N}}}{L^2(F)^d}^2  \le \int_{\Omega} s(c_0) \dx + C_f \Norm{\alpha}{L^{\infty}(\Omega)} |\Omega| T,
\end{split}
\end{equation*}
with no restrictions on the size of the time steps~$\{\tau_n\}_{n = 1}^{N}$.
\eremk
\end{remark}

\subsection{Matrix form}\label{SEC::MATRIX}

In order to write the BDF-LDG method~\eqref{EQ::VARIATIONAL} in matrix form, we define the following forms and functionals: for all~$\wh,\psih\in\Wh$ and~$\zh,\sigmah,\phih\in\Rh$,
\begin{alignat*}{2}
{\sf m}_h(\zh,\phih)&:=\left(\zh,\phih\right)_\Omega,\\
{\sf b}_h(\wh,\phih)&:=\left(\nabla_h \vh,  \phih \right)_\Omega-\sum_{F\in\Fho} \left(\jump{\wh}_{\sf N}, \mvl{\phih}_{1-\mvlparamF}\right)_F,\\
{\sf j}_h(\wh,\psih)&:=\sum_{F \in \Fho} \big(\mathsf{h}^{-1} \jump{\wh}_{\sf N}, \jump{\psih}_{\sf N} \big)_{F},\\
{\sf d}_h(\wh,\phih)&:=\sum_{K\in\Th}\left(\D\zh,\phih\right)_K,\\
{\sf n}_{h}(\wh;\sigmah,\phih)&:=\sum_{K\in\Th}\left(\D s''\left(u(\wh)\right)\sigmah, \phih\right)_K,\\
{\sf u}_{h}(\wh,\psih)&:=\left(u(\wh),\psih\right)_\Omega,\\
{\sf f}_{h}(\whp,\whq;\psih)&:=\left(f\left(u(\wh)\right),\psih\right)_\Omega.
\end{alignat*}
For fixed bases of~$\Wh$ and~$\Rh$, we denote by~$\cvW$, $\cvZ$, $\cvSigma$, and~$\cvR$ the coefficient vectors expressing~$\wh$, $\zh$, $\sigmah$, and~$\rh$, respectively, in terms of those bases. Additionally, we denote by~$M_I$, $B$, $J$, $M_D$, and~$A_{\LDG}$ the matrices associated with the bilinear forms~${\sf m}_h(\cdot,\cdot)$, ${\sf b}_h(\cdot,\cdot)$, ${\sf j}_h(\cdot,\cdot)$, ${\sf d}_h(\cdot,\cdot)$, and~$(\cdot,\cdot)_\LDG$, respectively, and by~${\mathcal N}$, ${\mathcal U}$, and~${\mathcal F}$ the operators associated with the non-linear functionals~${\sf n}_{h}(\cdot\,;\cdot,\cdot)$, ${\sf u}_{h}(\cdot,\cdot)$, and~${\sf f}_{h}(\cdot,\cdot\,;\cdot)$ respectively. Since the non-linear functional~${\sf n}_{h}(\cdot\,;\cdot,\cdot)$ is linear with respect to its second argument, we write~$\mathcal{N}(\cvW^{\npo}) \cvSigma^{\npo}$ instead of~$\mathcal{N}\big(\cvW^{\npo}; \cvSigma^{\npo}\big)$.

Then, the matrix form of the method in~\eqref{EQ::VARIATIONAL} is as follows:
\begin{subequations}
\label{EQ::MATRICIAL}
\begin{alignat}{3}
M_I \cvZ^{\npo} &= -B \cvW^{\npo}, \\
\mathcal{N}\big(\cvW^{\npo}\big) \cvSigma^{\npo}  & = M_D \cvZ^{\npo},\\
M_I \cvR^{\npo} & = M_D \cvSigma^{\npo},\\
\label{EQ::MATRICIAL-WH}
\varepsilon A_{\LDG} \cvW^{\npo} + \frac{1}{\tau_{n + 1} \beta}\,\mathcal{U}(\cvW^{\npo})  - B^T \cvR^{\npo} +  J \cvW^{(n + 1)} & = \mathcal{F}\big(\cvW^{\npo}) +\frac{1}{\tau_{n + 1}\beta}\,\sum_{j=1}^\bdfp a_j \mathcal{U}_{h}^{(n+1-j)},
\end{alignat}
\end{subequations}
where, in~\eqref{EQ::MATRICIAL-WH}, the terms~$\mathcal{U}_{h}^{(\kappa)}$ under the summation symbol represent the coefficient vector of~$\PiW p_0$ if~$\kappa = 0$, and
$\mathcal{U}_{h}^{(\kappa)} := \mathcal{U}\big(\cvW^{(\kappa)}\big)$ if~$\kappa > 0$.

As the mass matrix~$M_I$ is symmetric, positive definite, and block diagonal, it can be inverted efficiently and we can write 
\[
\cvZ^{\npo}=-M_I^{-1}B\cvW^{\npo}
\quad\text{and}\quad
\cvR^{\npo} = M_I^{-1}M_D \cvSigma^{\npo}.
\]
Then, system~\eqref{EQ::MATRICIAL} reduces to
\begin{subequations}
\label{EQ::MATRICIAL_COMPACT}
\begin{alignat}{3}
\label{EQ::MATRICIAL_COMPACT_SIGMA}
\mathcal{N}\big(\cvW^{\npo}\big) \cvSigma^{\npo}  & = -M_D M_I^{-1} B \cvW^{\npo},\\
\varepsilon A_{\LDG} \cvW^{\npo} + \frac{1}{\tau_{n + 1}\beta}\,\mathcal{U}(\cvW^{\npo}) 
- B^T M_I^{-1} M_D \cvSigma^{\npo} +  J \cvW^{(n + 1)} & = \mathcal{F} \big(\cvW^{\npo}\big) +\frac{1}{\tau_{n + 1} \beta}\,\sum_{j=1}^\bdfp a_j \mathcal{U}_{h}^{(n+1-j)}.
\end{alignat}
\end{subequations}

\begin{remark}[Reduced system with one unknown]\label{RMK::UNIQUESIGMA}
Given~$\cvW^{\npo}$, equation~\eqref{EQ::MATRICIAL_COMPACT_SIGMA} determines~$\cvSigma^{\npo}$ in a unique way. This follows from the uniform positive definiteness of~$\D s''(u(w))$; see~\eqref{EQ::UPD}. Then, using
$\cvSigma^{\npo}=-\big(\mathcal{N}\big(\cvW^{\npo}\big)\big)^{-1} M_D M_I^{-1}B \cvW^{\npo}$,
system~\eqref{EQ::MATRICIAL_COMPACT} can be reformulated  in the ~$\cvW^{\npo}$ unknown only as
\begin{equation}\label{EQ::MATRICIAL_W}
\begin{split}
\varepsilon A_{\LDG} \cvW^{\npo} + \frac{1}{\tau_{n + 1}\beta}\,\mathcal{U}(\cvW^{\npo})+B^TM_I^{-1}M_D \big(\mathcal{N}\big(\cvW^{\npo}\big)\big)^{-1} M_D M_I^{-1}B \cvW^{\npo}+J\cvW^{\npo}
\\
=\mathcal{F} \big(\cvW^{\npo}) +\frac{1}{\tau_{n + 1}\beta}\, \sum_{j=1}^\bdfp a_j \mathcal{U}_{h}^{(n+1-j)}.
\end{split}
\end{equation}
\eremk
\end{remark}

\subsection{Newton's iteration}\label{SEC::NEWTON}
To complete the presentation of the method, we derive explicitly the linear systems stemming from Newton's iteration. 
One of the most computationally expensive parts in each linear iteration is the evaluation of the non-linear terms in the multivariate function and its Jacobian matrix. However, since non-linearities in our method do not appear on interface integrals, these terms can be computed independently for each element~$K$ of~$\Th$, and their corresponding Jacobians are block diagonal. Consequently, our method has an inherent parallelizable structure.

We set, for convenience,~$C:=M_D M_I^{-1}B$ and
\[
\begin{split}
\GG_1\big(\cvSigma^{\npo},\cvW^{\npo}\big)&:=
\mathcal{N}\big(\cvW^{\npo}\big) \cvSigma^{\npo} +C \cvW^{\npo},\\[0.2cm]
\GG_2\big(\cvSigma^{\npo},\cvW^{\npo}\big)&:=
\varepsilon A_{\LDG} \cvW^{\npo} + \frac{1}{\tau_{n + 1}\beta}\,\mathcal{U}(\cvW^{\npo})  - C^T \cvSigma^{\npo} +  J \cvW^{(n + 1)}\\
& \qquad - \mathcal{F}\big(\cvW^{\npo}) -\frac{1}{\tau_{n + 1}\beta}\,\sum_{j=1}^{\bdfp}\,\mathcal{U}_{h}^{(n+1-j)}.
\end{split}
\]
Denote by $\diffsigma$\, and $\diffw$\, the differential operators with respect to~$\mathbf{\Sigma}$ and~$\mathbf{W}$, respectively. Omitting the temporal index~$n+1$,
the step~$k\to k+1$ of Newton's iteration applied to system~\eqref{EQ::MATRICIAL_COMPACT} reads as follows:
\[
\begin{split}
&\diffsigma\left(\GG_1\left(\cvSigma^{k},\cvW^{k}\right)\right)\big(\cvSigma^{k+1}-\cvSigma^{k}\big)
+\diffw\left(\GG_1\left(\cvSigma^{k},\cvW^{k}\right)\right)\big(\cvW^{k+1}-\cvW^{k}\big)
=-\GG_1\left(\cvSigma^{k},\cvW^{k}\right),
\\
&\diffsigma\left(\GG_2\left(\cvSigma^{k},\cvW^{k}\right)\right)\big(\cvSigma^{k+1}-\cvSigma^{k}\big)
+\diffw\left(\GG_2\left(\cvSigma^{k},\cvW^{k}\right)\right)\big(\cvW^{k+1}-\cvW^{k}\big)
=-\GG_2\left(\cvSigma^{k},\cvW^{k}\right).
\end{split}
\]
We compute
\[
\begin{split}
\diffsigma\left(\GG_1\left(\cvSigma^{k},\cvW^{k}\right)\right)
&=\mathcal{N}\big(\cvW^{k}\big),\\
\diffw\left(\GG_1\left(\cvSigma^{k},\cvW^{k}\right)\right)
&=\diffw\left(\mathcal{N}\big(\cvW^{k}\big)\right)\cvSigma^{k}+C,\\
\diffsigma\left(\GG_2\left(\cvSigma^{k},\cvW^{k}\right)\right)
&=-C^T,\\
\diffw\left(\GG_2\left(\cvSigma^{k},\cvW^{k}\right)\right)&=\varepsilon A_{\LDG}+\frac{1}{\tau_{n + 1}\beta}\,\diffw\left(\mathcal{U}(\cvW^{k})\right)+J
-\diffw\left(\mathcal{F} \big(\cvW^{k}\big)\right).
\end{split}
\]
Therefore, the Newton iteration applied to system~\eqref{EQ::MATRICIAL_COMPACT} is as follows:  
\begin{equation}\label{EQ::NEWTON}
\begin{split}
\mathcal{N}\big(\cvW^{k}\big) \cvSigma^{k+1}&=- \Big[\diffw\left(\mathcal{N}\big(\cvW^{k}\big)\right)\cvSigma^{k}+C \Big]\cvW^{k+1}+\left(\diffw\left(\mathcal{N}\big(\cvW^{k}\big)\right)\cvSigma^{k}\right)\cvW^{k},\\[0.5cm]
-C^T \cvSigma^{k+1}+\Big[&\varepsilon A_{\LDG}+\frac{1}{\tau_{n + 1}\beta}\,\diffw\left(\mathcal{U}(\cvW^{k})\right)+J-\diffw\left(\mathcal{F} \big(\cvW^{k}\big)\right)\Big]\cvW^{k+1}\\
&
=
\Big[\frac{1}{\tau_{n + 1}\beta}\,\diffw\left(\mathcal{U}(\cvW^{k})\right)-\diffw\left(\mathcal{F}\big(\cvW^{k}\big)\right)\Big]\cvW^{k} -\frac{1}{\tau_{n + 1} \beta}\, \Big[\mathcal{U}(\cvW^{k})-\sum_{j=1}^{\bdfp}\,\mathcal{U}_{h}^{(n+1-j)}\Big]
+ \mathcal{F} \big(\cvW^{k}\big).
\end{split}
\end{equation}

\begin{remark}[Newton's iteration for~\eqref{EQ::MATRICIAL_W}]
Using~\eqref{EQ::MATRICIAL_COMPACT_SIGMA}, one can eliminate~$\cvSigma^{k}$ from~\eqref{EQ::NEWTON} and obtain 
\begin{equation}\label{EQ::NEWTON_ONEFIELD}
\begin{split}
\Big[\varepsilon A_{\LDG}+ & \frac{1}{\tau_{n + 1}\beta}\,\diffw\left(\mathcal{U}(\cvW^{k})\right)
+C^T\left(\mathcal{N}\big(\cvW^{k}\big)\right)^{-1}C-\MM\cvW^{k}
+J-\diffw\left(\mathcal{F} \big(\cvW^{k}\big)\right) \Big]\cvW^{k+1} \\
\qquad
= &
\Big[\frac{1}{\tau_{n + 1}\beta}\,\diffw\left(\mathcal{U}(\cvW^{k})\right)-\diffw\left(\mathcal{F}\big(\cvW^{k}\big)\right)-\left(\MM\cvW^{k}\right)\Big]\cvW^{k} -\frac{1}{\tau_{n + 1}\beta}\,\Big[\mathcal{U}(\cvW^{k})-\sum_{j=1}^{\bdfp}\,\mathcal{U}_{h}^{(n+1-j)} \Big]
+ \mathcal{F} \big(\cvW^{k}\big),
\end{split}
\end{equation}
where~$\MM$ is the third-order tensor defined as
\[
\MM:=C^T\left(\mathcal{N}\big(\cvW^{k}\big)\right)^{-1}\diffw\left(\mathcal{N}\big(\cvW^{k}\big)\right)\left(\mathcal{N}\big(\cvW^{k}\big)\right)^{-1}C.
\]
The expression in~\eqref{EQ::NEWTON_ONEFIELD} can also be obtained by applying Newton's iteration directly to the reformulation of system~\eqref{EQ::MATRICIAL_COMPACT} given in~\eqref{EQ::MATRICIAL_W}.
\eremk
\end{remark}

\section{Numerical verification}\label{sec:verification}
In this section, we present some numerical tests to assess the accuracy of the proposed structure-preserving LDG method. We consider test cases in two space dimensions. In Section~\ref{sec:test_case_1}, we discuss the convergence properties of the scheme in space and time for a smooth exact solution. Then, in Section~\ref{sec:test_case_2}, we study the accuracy of the method in simulating a traveling-wave solution, providing a comparison with an existing interior penalty DG (IPDG) method that does not guarantee the positivity preservation~\cite{Corti_Bonizzoni_Dede_Quarteroni_Antonietti:2023}. 
All numerical simulations in this section are based on the \texttt{lymph} library \cite{antonietti_lymph_2024}, implementing 
DG methods on polytopic meshes. The polygonal meshes are constructed using PolyMesher~\cite{talischi_polymesher_2012}. All the constructed meshes satisfy the polytopic regularity assumptions of~\cite[Assumption 30 in~\S4.3]{Cangiani_Dong_Georgoulis:2017}. We use time meshes with uniform time step~$\tau$. 

Since the elements in the meshes we use in this section have a uniformly small number of edges, we neglect the dependence of the mesh-size function~$\mathsf{h}$ in~\eqref{EQN::DEF-h} on~$m_{K_*}$. We fix~$\theta = -1$ in~\eqref{EQN::DEF-h}, and set the stabilization parameter~$\eta_0 = 1$ in the definition of~$\eta_F$. For the Newton iterations, given a small tolerance~$\mathsf{tol}$, we adopt the following
stopping criterion: 
\begin{equation}
\label{eq:stopping-criteria}
\min\left\{\|w_h^{k+1}-w_h^k\|_{L^2(\Omega)}, |\mathsf{res}_{k+1}|\right\} \leq \mathsf{tol},
\end{equation} 
where $\mathsf{res}_{k+1}$ is the residual of the algebraic system~\eqref{EQ::MATRICIAL_W} for the approximation at the~$(k+1)$th Newton's iteration of~$w_h^{(n+1)}$.
In the convergence tests reported below, we measure the errors at the final time
\begin{equation*}
\Ec := \Norm{c(\cdot, T) - u(w_h^{(N)})}{L^2(\Omega)} \quad \text{ and } \quad \Esigma : = \Norm{\nabla c(\cdot, T) + \vsigma_h^{(N)}}{L^2(\Omega)^d}.
\end{equation*}

\subsection{Test case 1: Convergence analysis
}
\label{sec:test_case_1}
For the numerical tests in this section, 
we consider the space domain~$\Omega=(0,1)^2$
and homogeneous Neumann boundary conditions on~$\Gamma \times (0, T)$. 
For the non-linear Newton solver, we 
use the stopping criteria~\eqref{eq:stopping-criteria} with~$\mathsf{tol} = 10^{-16}$. 
The penalty parameter~$\varepsilon$ is set to~$0$.

\subsubsection*{Convergence properties of the space discretization}
\begin{figure}[t!]
    \begin{subfigure}[b]{0.25\textwidth}
    \includegraphics[width=\textwidth]{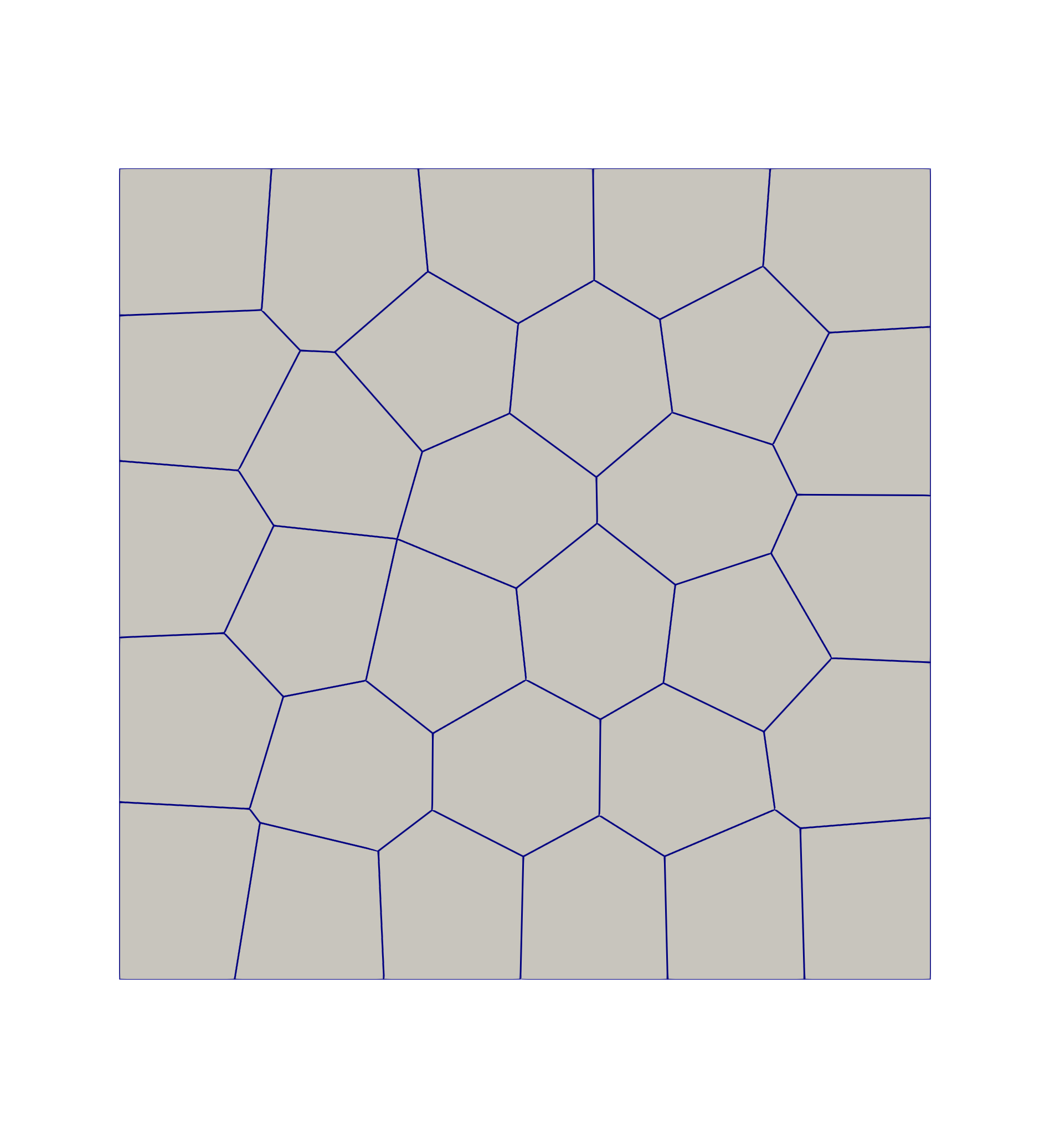}
    \caption{Mesh of 30 elements.}
    \end{subfigure}%
    \begin{subfigure}[b]{0.25\textwidth}
    \includegraphics[width=\textwidth]{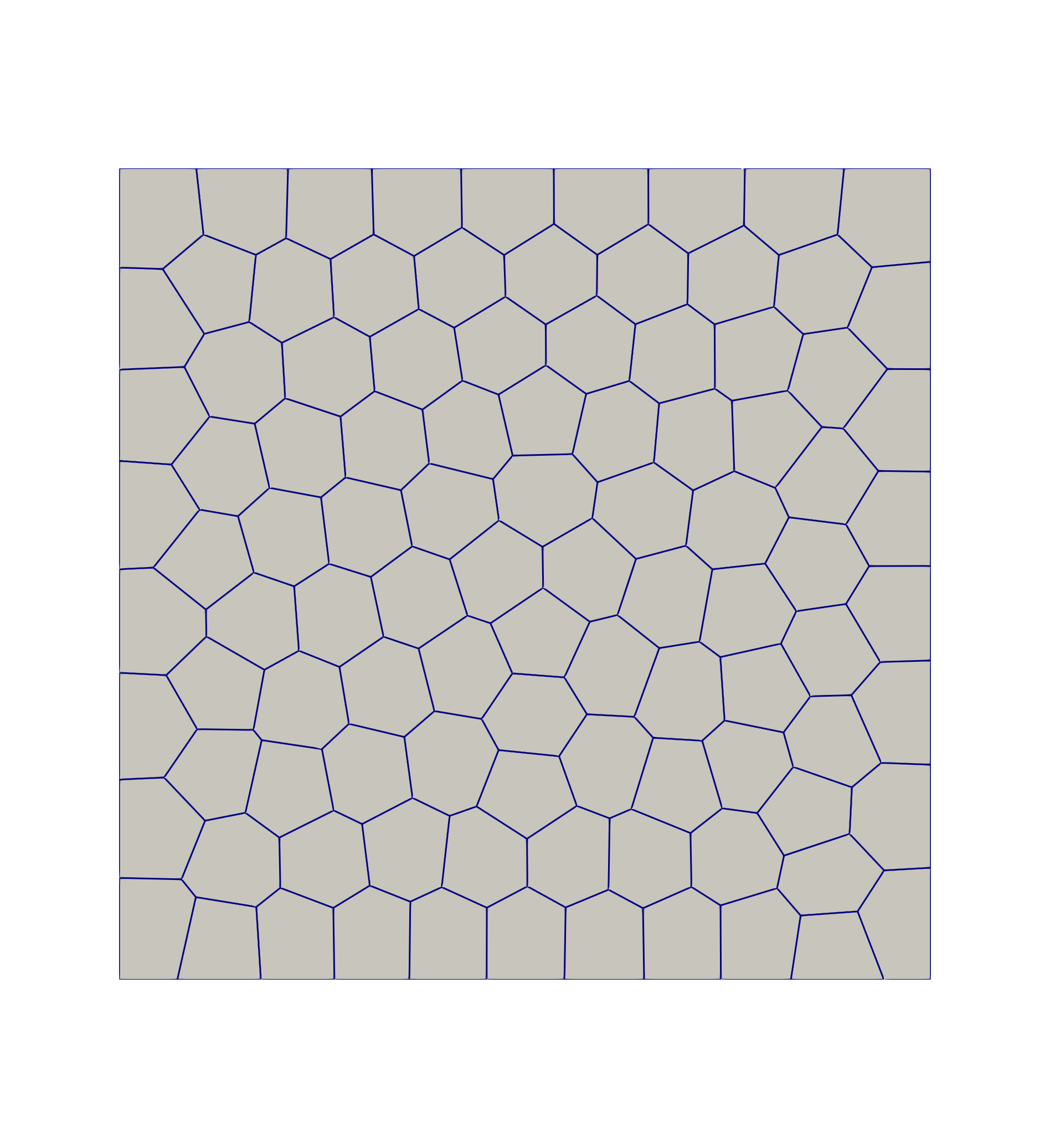}
    \caption{Mesh of 100 elements.}
    \end{subfigure}%
    \begin{subfigure}[b]{0.25\textwidth}
    \includegraphics[width=\textwidth]{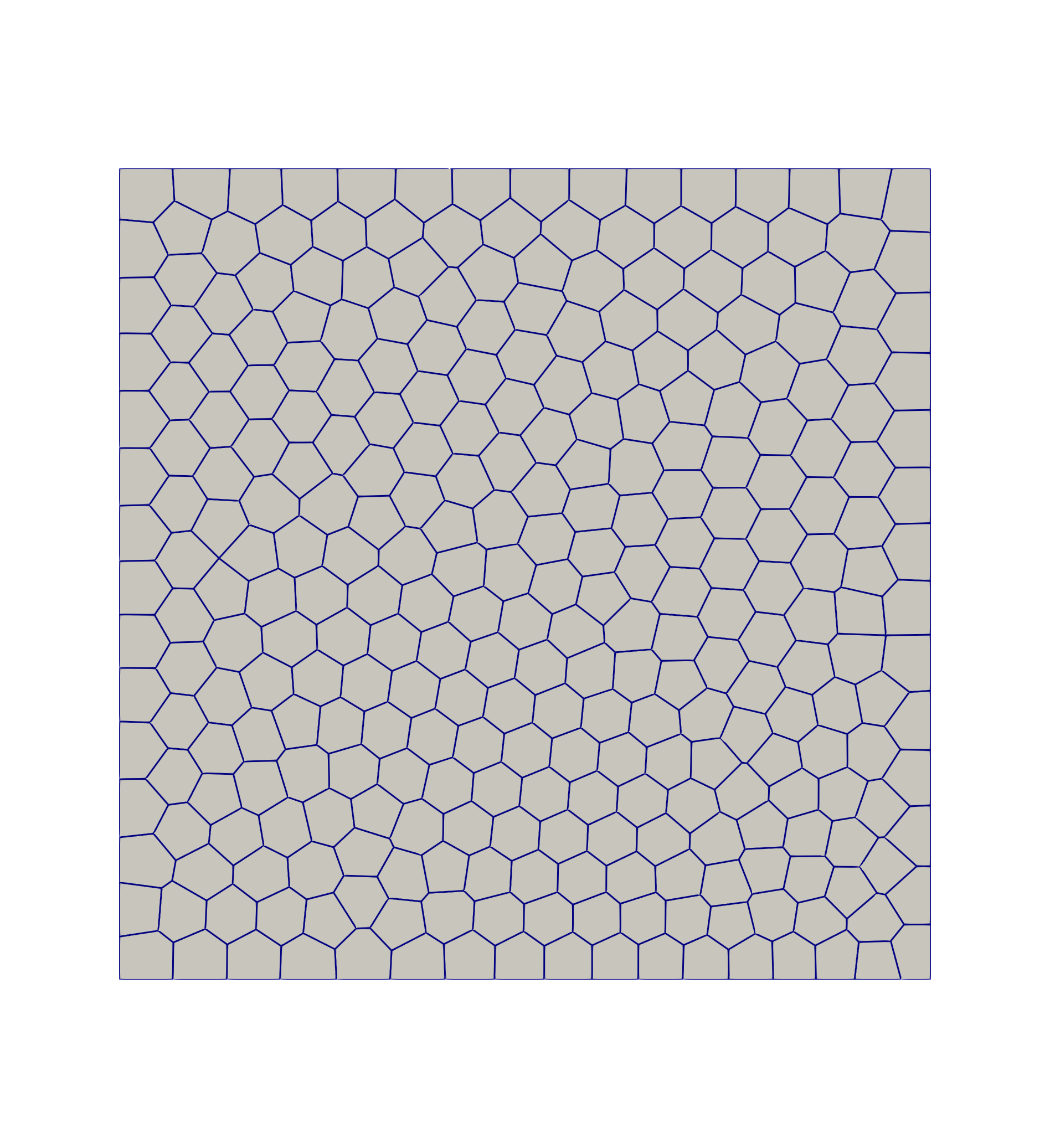}
    \caption{Mesh of 300 elements.}
    \end{subfigure}%
    \begin{subfigure}[b]{0.25\textwidth}
    \includegraphics[width=\textwidth]{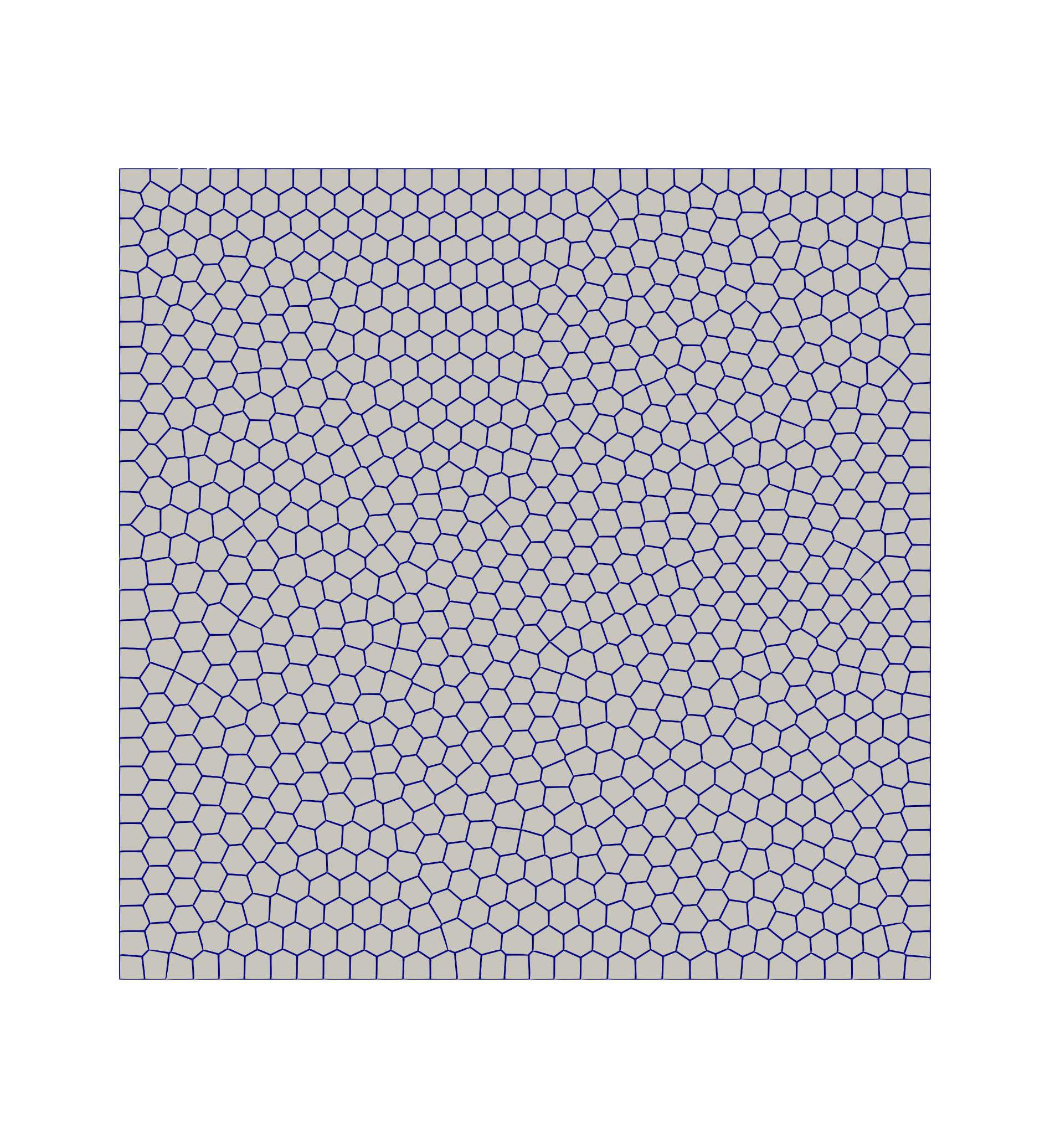}
    \caption{Mesh of 1000 elements.}
    \end{subfigure}%
    
    \caption{Test case 1: meshes with different levels of refinement used in the $h$-convergence test.}
    \label{fig:meshes}
\end{figure}
For this test, we use the BDF1 time-stepping scheme (implicit Euler) for the time discretization. We select the diffusion tensor~$\mathbf{D}=\mathbb{I}_2$, where~$\mathbb{I}_2$ represents the identity matrix of size~$2$, and the reaction coefficient~$\alpha=1$, and consider the FK equation with an additional source term on the right-hand side so that problem~\eqref{EQN::FISHER-KPP} admits
the following manufactured exact solution:
\begin{equation*}
    c(x,y,t) = \frac{1}{4}\left(\cos(2 \pi x)\cos(2 \pi y)+2\right) (1-t).
\end{equation*}
The choice of a solution that is linear in time allows us to neglect the error due to the BDF scheme, 
highlighting the properties of the space discretization. 
\par
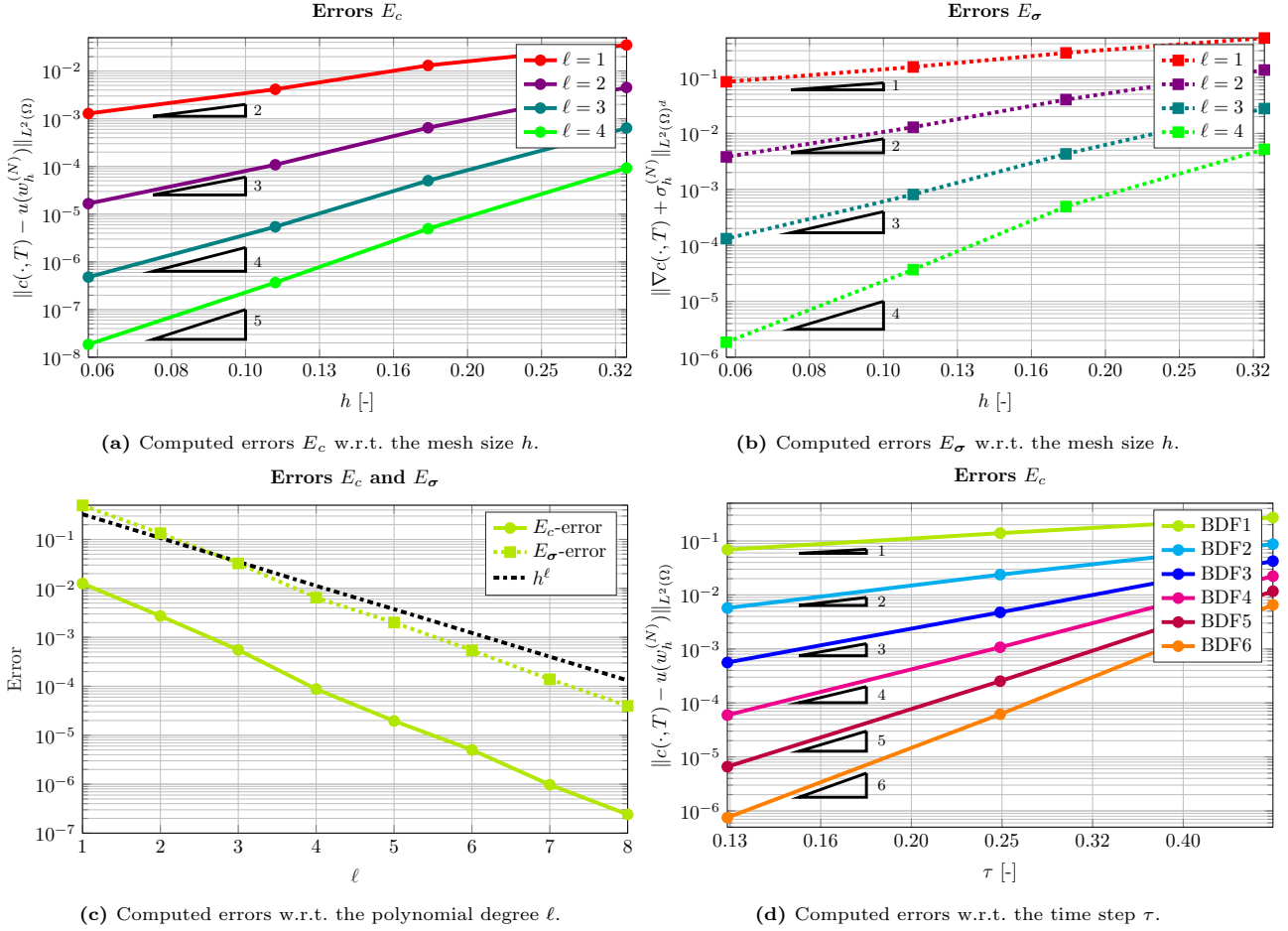
\begin{figure}[t!]
    \begin{subfigure}[b]{0.5\textwidth}
          \resizebox{\textwidth}{!}{\definecolor{mycolor2}{rgb}{0.00000,1.00000,1.00000}%
\begin{tikzpicture}
\begin{axis}[%
width=3.875in,
height=2.30in,
at={(2.6in,1.099in)},
scale only axis,
xmode=log,
xmin=0.061309,
xmax=3.2726e-1,
xminorticks=true,
xlabel = {$h$ [-]},
ylabel = {$\|c(\cdot,T)-u(w_h^{(N)})\|_{L^2(\Omega)}$},
xticklabel={\pgfmathparse{exp(\tick)}\pgfmathprintnumber{\pgfmathresult}},
x tick label style={
/pgf/number format/.cd, fixed, fixed zerofill,
precision=2},
ymode=log,
ymin=1e-8,
ymax=0.05,
yminorticks=true,
axis background/.style={fill=white},
title style={font=\bfseries},
title={Errors $E_{c}$},
xmajorgrids,
xminorgrids,
ymajorgrids,
yminorgrids,
legend style={legend cell align=left, align=left, draw=white!15!black}
]
              
\addplot [color=red, line width=2.0pt, mark=*]
  table[row sep=crcr]{%
3.2726e-1  3.5231e-2\\
1.7647e-1  1.3118e-2\\
1.0969e-1  4.1462e-3\\
6.1309e-2  1.2862e-3\\
};
\addlegendentry{$\ell=1$}

\addplot [color=violet, line width=2.0pt, mark=*]
  table[row sep=crcr]{%
3.2726e-1  4.5231e-3\\
1.7647e-1  6.4734e-4\\
1.0969e-1  1.0854e-4\\
6.1309e-2  1.6620e-5\\
};
\addlegendentry{$\ell=2$}

\addplot [color=teal, line width=2.0pt, mark=*]
  table[row sep=crcr]{%
3.2726e-1  6.3913e-4\\
1.7647e-1  5.0416e-5\\
1.0969e-1  5.4358e-6\\
6.1309e-2  4.7878e-7\\
};
\addlegendentry{$\ell=3$}

\addplot [color=green, line width=2.0pt, mark=*]
  table[row sep=crcr]{%
3.2726e-1  9.2535e-5\\
1.7647e-1  4.9388e-6\\
1.0969e-1  3.6809e-7\\
6.1309e-2  1.8576e-8\\
};
\addlegendentry{$\ell=4$}



\node[right, align=left, text=black, font=\footnotesize]
at (axis cs:0.1005,0.0015) {$2$};

\addplot [color=black, line width=1.5pt]
  table[row sep=crcr]{%
0.100   0.002\\
0.075   0.001125\\
0.100   0.001125\\
0.100   0.002\\
};

\node[right, align=left, text=black, font=\footnotesize]
at (axis cs:0.1005,4e-5) {$3$};

\addplot [color=black, line width=1.5pt]
  table[row sep=crcr]{%
0.100   6e-05\\
0.075   2.53e-05\\
0.100   2.53e-05\\
0.100   6e-05\\
};

\node[right, align=left, text=black, font=\footnotesize]
at (axis cs:0.1005,1e-6) {$4$};

\addplot [color=black, line width=1.5pt]
  table[row sep=crcr]{%
0.100   2e-06\\
0.075   6.328e-7\\
0.100   6.328e-7\\
0.100   2e-6\\
};

\node[right, align=left, text=black, font=\footnotesize]
at (axis cs:0.1005,6e-8) {$5$};

\addplot [color=black, line width=1.5pt]
  table[row sep=crcr]{%
0.100   1e-7\\
0.075   2.37e-8\\
0.100   2.37e-8\\
0.100   1e-7\\
};



\end{axis}
\end{tikzpicture}%
}
          \caption{Computed errors $E_{c}$ w.r.t.~the mesh size~$h$.}
        \label{fig:errors2D_h_L2}
    \end{subfigure}%
    \begin{subfigure}[b]{0.5\textwidth}
        \resizebox{\textwidth}{!}{\definecolor{mycolor2}{rgb}{0.00000,1.00000,1.00000}%
\begin{tikzpicture}

\begin{axis}[%
width=3.875in,
height=2.30in,
at={(2.6in,1.099in)},
scale only axis,
xmode=log,
xmin=6.1309e-2,
xmax=3.2726e-1,
xminorticks=true,
xlabel = {$h$ [-]},
ylabel = {$\|\nabla c(\cdot,T) + \sigma_h^{(N)}\|_{L^2(\Omega)^d}$},
xticklabel={\pgfmathparse{exp(\tick)}\pgfmathprintnumber{\pgfmathresult}},
x tick label style={
/pgf/number format/.cd, fixed, fixed zerofill,
precision=2},
ymode=log,
ymin=1e-6,
ymax=0.51,
yminorticks=true,
axis background/.style={fill=white},
title style={font=\bfseries},
title={Errors $E_{\boldsymbol{\sigma}}$},
xmajorgrids,
xminorgrids,
ymajorgrids,
yminorgrids,
legend style={legend cell align=left, align=left, draw=white!15!black}
]
              
\addplot [color=red, line width=2.0pt, mark=square*, dotted, mark options=solid]
  table[row sep=crcr]{%
3.2726e-1  5.0613e-1\\
1.7647e-1  2.7437e-1\\
1.0969e-1  1.5371e-1\\
6.1309e-2  8.3484e-2\\
};
\addlegendentry{$\ell=1$}

\addplot [color=violet, line width=2.0pt, mark=square*, dotted, mark options=solid]
  table[row sep=crcr]{%
3.2726e-1  1.3560e-1\\
1.7647e-1  4.0152e-2\\
1.0969e-1  1.2905e-2\\
6.1309e-2  3.8056e-3\\
};
\addlegendentry{$\ell=2$}

\addplot [color=teal, line width=2.0pt, mark=square*, dotted, mark options=solid]
  table[row sep=crcr]{%
3.2726e-1  2.7800e-2\\
1.7647e-1  4.3081e-3\\
1.0969e-1  8.1153e-4\\
6.1309e-2  1.3152e-4\\
};
\addlegendentry{$\ell=3$}

\addplot [color=green, line width=2.0pt, mark=square*, dotted, mark options=solid]
  table[row sep=crcr]{%
3.2726e-1  5.2000e-3\\
1.7647e-1  4.9388e-4\\
1.0969e-1  3.6809e-5\\
6.1309e-2  1.8576e-6\\
};
\addlegendentry{$\ell=4$}



\node[right, align=left, text=black, font=\footnotesize]
at (axis cs:0.1005,0.072) {$1$};

\addplot [color=black, line width=1.5pt]
  table[row sep=crcr]{%
0.100   0.08\\
0.075   0.06\\
0.100   0.06\\
0.100   0.08\\
};

\node[right, align=left, text=black, font=\footnotesize]
at (axis cs:0.1005,0.006) {$2$};

\addplot [color=black, line width=1.5pt]
  table[row sep=crcr]{%
0.100   0.008\\
0.075   0.0045\\
0.100   0.0045\\
0.100   0.008\\
};

\node[right, align=left, text=black, font=\footnotesize]
at (axis cs:0.1005,2.5e-4) {$3$};

\addplot [color=black, line width=1.5pt]
  table[row sep=crcr]{%
0.100   4e-04\\
0.075   1.68e-04\\
0.100   1.68e-04\\
0.100   4e-04\\
};

\node[right, align=left, text=black, font=\footnotesize]
at (axis cs:0.1005,6e-6) {$4$};

\addplot [color=black, line width=1.5pt]
  table[row sep=crcr]{%
0.100   1e-05\\
0.075   3.164e-06\\
0.100   3.164e-06\\
0.100   1e-05\\
};





\end{axis}
\end{tikzpicture}%
}
          \caption{Computed errors $E_{\boldsymbol{\sigma}}$
          w.r.t.~the mesh size $h$.}
        \label{fig:errors2D_h_dG}
    \end{subfigure}%
    \\[0.5mm]
    \begin{subfigure}[b]{0.5\textwidth}
          \resizebox{\textwidth}{!}{\definecolor{mylime}{rgb}{0.70000,0.90000,0.00000}%

\begin{tikzpicture}
\begin{axis}[%
width=3.875in,
height=2.33in,
at={(2.6in,1.099in)},
scale only axis,
xmin=1,
xmax=8,
xlabel style={font=\color{white!15!black}},
xlabel={$\ell$},
ymode=log,
ymin=1e-7,
ymax=5e-1,
yminorticks=true,
ylabel style={font=\color{white!15!black}},
ylabel={Error},
axis background/.style={fill=white},
title style={font=\bfseries},
title={Errors $E_{c}$ and $E_{\boldsymbol{\sigma}}$},
xmajorgrids,
xminorgrids,
ymajorgrids,
yminorgrids,
legend style={legend cell align=left, align=left, draw=white!15!black}
]

\addplot [color=mylime, line width=2.0pt, mark=*]
  table[row sep=crcr]{%
    1   1.26e-02 \\
    2   2.73e-03 \\
    3   5.54e-04 \\
    4   8.76e-05 \\
    5   1.96e-05 \\
    6   5.02e-06 \\
    7   9.74e-07 \\
    8   2.41e-07 \\
};
\addlegendentry{$E_{c}$-error}

\addplot [color=mylime, line width=2.0pt, mark=square*, dotted, mark options=solid]
  table[row sep=crcr]{%
    1   4.95e-01 \\
    2   1.34e-01 \\
    3   3.23e-02 \\
    4   6.46e-03 \\
    5   2.01e-03 \\
    6   5.42e-04 \\
    7   1.39e-04 \\
    8   3.96e-05 \\
};
\addlegendentry{$E_{\boldsymbol{\sigma}}$-error}

\addplot [color=black, line width=2.0pt, dashdotted, mark options=solid]
  table[row sep=crcr]{%
    1   3.2726e-01 \\
    2   1.0709e-01 \\
    3   3.5049e-02 \\
    4   1.1270e-02 \\
    5  	3.7537e-03 \\
    6  	1.2284e-03 \\
    7  	4.0202e-04 \\
    8  	1.3156e-04 \\
};
\addlegendentry{$h^\ell$}

\end{axis}
\end{tikzpicture}%
}
          \caption{Computed errors 
          w.r.t.~the polynomial degree~$\ell$.}
        \label{fig:errors2D_p}
    \end{subfigure}%
    \begin{subfigure}[b]{0.5\textwidth}
        \resizebox{\textwidth}{!}{\definecolor{mylime}{rgb}{0.70000,0.90000,0.00000}%

\begin{tikzpicture}

\begin{axis}[%
width=3.875in,
height=2.30in,
at={(2.6in,1.099in)},
scale only axis,
xmode=log,
xmin=0.125,
xmax=0.5,
xminorticks=true,
xlabel = {$\tau$ [-]},
ylabel = {$\|c(\cdot,T)-u(w_h^{(N)})\|_{L^2(\Omega)}$},
xticklabel={\pgfmathparse{exp(\tick)}\pgfmathprintnumber{\pgfmathresult}},
x tick label style={
/pgf/number format/.cd, fixed, fixed zerofill,
precision=2},
ymode=log,
ymin=5e-7,
ymax=5e-1,
yminorticks=true,
axis background/.style={fill=white},
title style={font=\bfseries},
title={Errors $E_{c}$},
xmajorgrids,
xminorgrids,
ymajorgrids,
yminorgrids,
legend style={legend cell align=left, align=left, draw=white!15!black}
]
              
\addplot [color=mylime, line width=2.0pt, mark=*]
  table[row sep=crcr]{%
5.0000e-1  2.7157e-01\\
2.5000e-1  1.3923e-01\\
1.2500e-1  6.9022e-02\\
};
\addlegendentry{BDF1}

\addplot [color=cyan, line width=2.0pt, mark=*]
  table[row sep=crcr]{%
5.0000e-1  8.7157e-02\\
2.5000e-1  2.3742e-02\\
1.2500e-1  5.7182e-03\\
};
\addlegendentry{BDF2}

\addplot [color=blue, line width=2.0pt, mark=*]
  table[row sep=crcr]{%
5.0000e-1  4.2298e-02\\
2.5000e-1  4.7559e-03\\
1.2500e-1  5.6078e-04\\
};
\addlegendentry{BDF3}

\addplot [color=magenta, line width=2.0pt, mark=*]
  table[row sep=crcr]{%
5.0000e-1  2.2286e-02\\
2.5000e-1  1.0725e-03\\
1.2500e-1  5.9493e-05\\
};
\addlegendentry{BDF4}

\addplot [color=purple, line width=2.0pt, mark=*]
  table[row sep=crcr]{%
5.0000e-1  1.1700e-02\\
2.5000e-1  2.5245e-04\\
1.2500e-1  6.5816e-06\\
};
\addlegendentry{BDF5}

\addplot [color=orange, line width=2.0pt, mark=*]
  table[row sep=crcr]{%
5.0000e-1  6.5496e-03\\
2.5000e-1  6.1654e-05\\
1.2500e-1  7.5060e-07\\
};
\addlegendentry{BDF6}

\node[right, align=left, text=black, font=\footnotesize]
at (axis cs:0.180,0.065) {$1$};

\addplot [color=black, line width=1.5pt]
  table[row sep=crcr]{%
0.178   7.0000e-02\\
0.150   5.8989e-02\\
0.178   5.8989e-02\\
0.178   7.0000e-02\\
};

\node[right, align=left, text=black, font=\footnotesize]
at (axis cs:0.1800,0.0075) {$2$};

\addplot [color=black, line width=1.5pt]
  table[row sep=crcr]{%
0.178   9.0000e-03\\
0.150   6.3912e-03\\
0.178   6.3912e-03\\
0.178   9.0000e-03\\
};

\node[right, align=left, text=black, font=\footnotesize]
at (axis cs:0.18,1e-3) {$3$};

\addplot [color=black, line width=1.5pt]
  table[row sep=crcr]{%
0.178   1.2500e-03\\
0.150   7.4804e-04\\
0.178   7.4804e-04\\
0.178   1.2500e-03\\
};

\node[right, align=left, text=black, font=\footnotesize]
at (axis cs:0.1800,1.4e-4) {$4$};

\addplot [color=black, line width=1.5pt]
  table[row sep=crcr]{%
0.178   2.0000e-04\\
0.150   1.0086e-04\\
0.178   1.0086e-04\\
0.178   2.0000e-04\\
};

\node[right, align=left, text=black, font=\footnotesize]
at (axis cs:0.1800,2e-5) {$5$};

\addplot [color=black, line width=1.5pt]
  table[row sep=crcr]{%
0.178   3.0000e-05\\
0.150   1.2749e-05\\
0.178   1.2749e-05\\
0.178   3.0000e-05\\
};

\node[right, align=left, text=black, font=\footnotesize]
at (axis cs:0.1800,3e-6) {$6$};

\addplot [color=black, line width=1.5pt]
  table[row sep=crcr]{%
0.178   5.0000e-06\\
0.150   1.7906e-06\\
0.178   1.7906e-06\\
0.178   5.0000e-06\\
};

\end{axis}
\end{tikzpicture}%
}
          \caption{Computed errors 
          w.r.t.~the time step~$\tau$.}
        \label{fig:errors2D_dt}
    \end{subfigure}%
    \caption{Test case 1: computed errors and convergence rates. }
    \label{fig:errors2D}
\end{figure}

We perform a convergence test keeping fixed the polynomial degree of the space approximation $\ell=1,2,3,4$ and using, for each degree, different mesh refinements with number of elements $(N_\mathrm{el}= 30,100,300,1000)$. The polygonal meshes used in this test case are reported in Figure \ref{fig:meshes}. We take~$\tau= 10^{-3}$ and a final time $T=5\times10^{-2}$. In Figures \ref{fig:errors2D_h_L2} and \ref{fig:errors2D_h_dG}, we report the computed errors~$\Ec$ and~$\Esigma$. 
We observe a decrease in the error with quasi-optimal 
convergence,
namely, of order~$\mathcal{O}(h^\ell)$ 
for~$\Esigma$, and of order~$\mathcal{O}(h^{\ell+1})$ for~$\Ec$. 
\par
Then, we develop a convergence analysis with respect to the polynomial degree $\ell$. In this case, we consider a coarse mesh with~$30$ polygonal elements, 
and we take $\tau=10^{-5}$ and $T=2.5\times10^{-4}$. 
The errors~$\Ec$ and~$\Esigma$
are reported in Figure~\ref{fig:errors2D_p}, where we also report the line of convergence rate $h^\ell$. Comparing the results, we can observe a spectral convergence in the polynomial degree $\ell$.
%
\subsubsection*{Convergence properties of the time discretization}
We also test the convergence in time of the fully discrete scheme by comparing the six stable BDF schemes. For this test, we select the diffusion tensor~$\mathbf{D}=10^{-3}\mathbb{I}_2$ and the reaction coefficient~$\alpha=1$, and we slightly modify the time dependence in the manufactured exact solution with respect to the previous experiment to highlight the properties of the time discretization: 
\begin{equation*}
    c(x,y,t) = \frac{1}{4}\left(\cos(2 \pi x)\cos(2 \pi y)+2\right) e^{-t}.
\end{equation*}
We use a fixed mesh of~$300$ elements ($h \approx 0.1097$) and a polynomial degree~$\ell=4$ for the space discretization. This guarantees that the errors induced by the space discretization are sufficiently small to highlight the convergence properties of the time integration methods. 
\par
The convergence test uses three different time steps $\tau = 0.500,\, 0.250,\,0.125$, and a final time $T=2$.
In Figure \ref{fig:errors2D_dt}, we report the computed errors~$\Ec$.
Coherently with the theoretical results on BDF methods, we observe a decrease in the error of order~$\mathcal{O}(\tau^\bdfp)$ for the general \BDF{\bdfp} method.
\subsection{Test case 2: Traveling-wave solution}
\label{sec:test_case_2}
In this section, we emphasize the critical role of
the strong positivity-preserving property of the 
BDF-LDG formulation~\eqref{EQ::MATRICIAL_COMPACT} in accurately simulating
a traveling-wave solution. We fix an anisotropic diffusion tensor $\mathbf{D}=\dext\mathbb{I}_2$ and a constant reaction coefficient~$\alpha$, and we consider a solution of the form
\begin{equation*} 
    c(x,y,t) = \psi(x-vt) = \psi(\xi),
\end{equation*}
where~$v$ is wave speed depending on~$\dext$ and~$\alpha$ defined by~$v := 5\sqrt{{\alpha \dext}/{6}}$.
Substituting~$c$ 
in the FK equation, we 
obtain the following equivalent system of ordinary differential equations:
\begin{equation}
\label{eq:wave_eq}
    \begin{dcases}
        \chi'(\xi) = -\dfrac{v}{\dext}\chi(\xi) + \dfrac{\alpha}{\dext}\psi(\xi)(\psi(\xi)-1) & \xi\in(0,T), \\
        \psi'(\xi) = \chi(\xi) & \xi\in(0,T).
        \\
    \end{dcases}
\end{equation}
The analytical solution of problem \eqref{eq:wave_eq} is given by (see~\cite[\S7.2]{wen-shan_exact_2006})
\begin{equation*}
    \psi(\xi) = \dfrac{1}{4}\left(1+\tanh\left(8-\sqrt{\dfrac{\alpha}{24\dext}}\xi\right)\right)^2.
\end{equation*}
This solution 
satisfies a homogeneous Neumann boundary condition at the limits~$\xi\rightarrow\pm\infty$, which is equivalent to~$x\rightarrow\pm\infty$ for each fixed value of~$t\in(0,T)$. 
The homogeneous Neumann boundary condition is also respected in the $y$-direction, as the exact solution~$c$ is independent of~$y$. In this simulation, we consider a rectangular space domain~$\Omega = (0,3)\times(0,1)$, and the final time~$T=10$. We impose homogeneous Neumann boundary conditions not only at~$y=0$ and~$y=1$, but also at~$x=0$ and~$x=3$.
\par
Concerning the physical parameters, we fix $\dext=10^{-3}$ and~$\alpha=1$, 
with associated velocity 
$v\simeq 6.45\times 10^{-2}$. 
For the non-linear Newton solver, we 
use the stopping criteria~\eqref{eq:stopping-criteria} with~$\mathsf{tol} = 10^{-10}$.  

\subsubsection*{Impact of the penalty 
parameter~$\varepsilon$}
In the first test, we analyze the impact of the 
penalty parameter~$\varepsilon$ on the error convergence. To test this, we fix the polynomial degree~$\ell=5$, and the mesh
consisting of 50 elements generated with PolyMesher~\cite{talischi_polymesher_2012}. 
We consider as time step~$\tau = 2.5\times 10^{-2}$, and as final time $T=10$. In this case, we test two different choices of \BDF{\bdfp} scheme, namely, $\bdfp = 1$ and~$\bdfp = 6$.
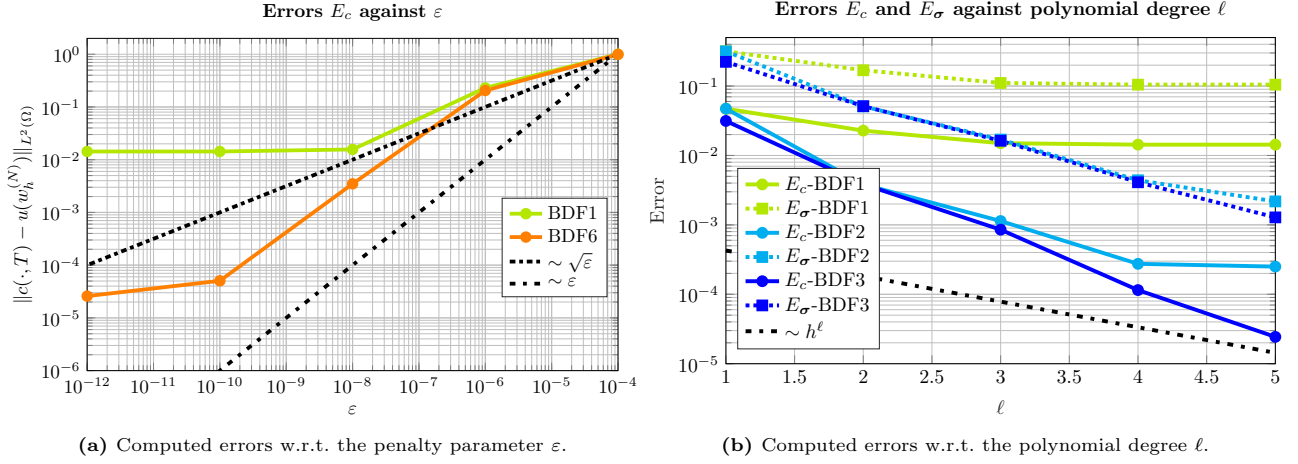
\begin{figure}[t!]
    \begin{subfigure}[b]{0.5\textwidth}
        \resizebox{\textwidth}{!}{\definecolor{mylime}{rgb}{0.70000,0.90000,0.00000}%

\begin{tikzpicture}

\begin{axis}[%
width=3.875in,
height=2.425in,
at={(2.6in,1.099in)},
scale only axis,
xmode=log,
xmin=1e-12,
xmax=1e-4,
xminorticks=true,
xlabel = {$\varepsilon$},
ylabel = {$\| c(\cdot,T) - u(w_h^{(N)})\|_{L^2(\Omega)}$},
ymode=log,
ymin=1e-6,
ymax=2,
yminorticks=true,
axis background/.style={fill=white},
title style={font=\bfseries},
title={Errors $E_c$ against $\varepsilon$},
xmajorgrids,
xminorgrids,
ymajorgrids,
yminorgrids,
legend style={at={(0.98,0.52)},legend cell align=left, draw=white!15!black}
]

\addplot [color=mylime, line width=2.0pt, mark=*]
  table[row sep=crcr]{%
    1.0000e-04   9.9827e-01 \\
    1.0000e-06   2.2974e-01 \\
    1.0000e-08   1.5645e-02 \\
    1.0000e-10   1.4315e-02 \\
    1.0000e-12   1.4315e-02 \\
};
\addlegendentry{BDF1}

\addplot [color=orange, line width=2.0pt, mark=*]
  table[row sep=crcr]{%
    1.0000e-04   9.9457e-01 \\
    1.0000e-06   2.0506e-01 \\
    1.0000e-08   3.4858e-03 \\
    1.0000e-10   5.0530e-05 \\
    1.0000e-12   2.5888e-05 \\
};
\addlegendentry{BDF6}

\addplot [color=black, line width=2.0pt, dashdotted, mark options=solid]
  table[row sep=crcr]{%
    1.0000e-04   1.0000e-00 \\
    1.0000e-06   1.0000e-01 \\
    1.0000e-08   1.0000e-02 \\
    1.0000e-10   1.0000e-03 \\
    1.0000e-12   1.0000e-04 \\
};
\addlegendentry{$\sim \sqrt{\varepsilon}$}

\addplot [color=black, line width=2.0pt, loosely dashdotted, mark options=solid]
  table[row sep=crcr]{%
    1.0000e-04   1.00e-00 \\
    1.0000e-06   1.00e-02 \\
    1.0000e-08   1.00e-04 \\
    1.0000e-10   1.00e-06 \\
    1.0000e-12   1.00e-08 \\
};
\addlegendentry{$\sim \varepsilon$}

\end{axis}
\end{tikzpicture}%
}
        \caption{Computed errors 
        w.r.t.~the penalty parameter~$\varepsilon$.}
        \label{fig:errorswaves_eps}
    \end{subfigure}%
    \begin{subfigure}[b]{0.5\textwidth}
        \resizebox{\textwidth}{!}{\definecolor{mylime}{rgb}{0.70000,0.90000,0.00000}%

\begin{tikzpicture}
\begin{axis}[%
width=3.875in,
height=2.30in,
at={(2.6in,1.099in)},
scale only axis,
xmin=1,
xmax=5,
xlabel style={font=\color{white!15!black}},
xlabel={$\ell$},
ymode=log,
ymin=1e-5,
ymax=5e-1,
yminorticks=true,
ylabel style={font=\color{white!15!black}},
ylabel={Error},
axis background/.style={fill=white},
title style={font=\bfseries},
title={Errors $E_c$ and $E_{\boldsymbol{\sigma}}$ against polynomial degree $\ell$},
xmajorgrids,
xminorgrids,
ymajorgrids,
yminorgrids,
    legend style={at={(0.28,0.60)},legend cell align=left, draw=white!15!black}
]

\addplot [color=mylime, line width=2.0pt, mark=*]
  table[row sep=crcr]{%
    1   4.7244e-02 \\
    2   2.2800e-02 \\
    3   1.5033e-02 \\
    4   1.4315e-02 \\
    5   1.4315e-02 \\
};
\addlegendentry{$E_c$-BDF1}

\addplot [color=mylime, line width=2.0pt, mark=square*, dotted, mark options=solid]
  table[row sep=crcr]{%
    1   3.1695e-01 \\
    2   1.6902e-01 \\
    3   1.1085e-01 \\
    4   1.0482e-01 \\
    5   1.0479e-01 \\
};
\addlegendentry{$E_{\boldsymbol{\sigma}}$-BDF1}

\addplot [color=cyan, line width=2.0pt, mark=*]
  table[row sep=crcr]{%
    1   4.7244e-02 \\
    2   3.8302e-03 \\
    3   1.1364e-03 \\
    4   2.7460e-04 \\
    5   2.5060e-04 \\
};
\addlegendentry{$E_c$-BDF2}

\addplot [color=cyan, line width=2.0pt, mark=square*, dotted, mark options=solid]
  table[row sep=crcr]{%
    1   3.1695e-01 \\
    2   5.0956e-02 \\
    3   1.6968e-02 \\
    4   4.3919e-03 \\
    5   2.1783e-03 \\
};
\addlegendentry{$E_{\boldsymbol{\sigma}}$-BDF2}

\addplot [color=blue, line width=2.0pt, mark=*]
  table[row sep=crcr]{%
    1   3.1441e-02 \\
    2   3.8553e-03 \\
    3   8.5273e-04 \\
    4   1.1501e-04 \\
    5   2.4336e-05 \\
};
\addlegendentry{$E_c$-BDF3}

\addplot [color=blue, line width=2.0pt, mark=square*, dotted, mark options=solid]
  table[row sep=crcr]{%
    1   2.2338e-01 \\
    2   5.1275e-02 \\
    3   1.6290e-02 \\
    4   4.1105e-03 \\
    5  	1.2826e-03 \\
};
\addlegendentry{$E_{\boldsymbol{\sigma}}$-BDF3}

\addplot [color=black, line width=2.0pt, loosely dashdotted, mark options=solid]
  table[row sep=crcr]{%
    1   4.2784e-4 \\
    2   1.8304e-4 \\
    3   7.8314e-5 \\
    4   3.3506e-5 \\
    5  	1.4335e-5 \\
};
\addlegendentry{$\sim h^\ell$}




\end{axis}
\end{tikzpicture}%
}
        \caption{Computed errors 
        w.r.t.~the polynomial degree~$\ell$.}
        \label{fig:errorswaves_ell}
    \end{subfigure}%
    \caption{Test case 2: computed errors 
    w.r.t.~the penalty parameter~$\varepsilon$ (a), and w.r.t.~the polynomial degree~$\ell$ (b).}
\end{figure}
\par
As we can observe in Figure \ref{fig:errorswaves_eps}, increasing the
penalty parameter introduces an error that seems to be of order~$\mathcal{O}(\varepsilon)$. 
Such a behavior is especially evident in the case of the \BDF{6} discretization. For the \BDF{1} scheme, a plateau is observed due to the dominance of the time discretization error.
In~\cite[\S5.1]{Gomez_Jungel_Perugia:2024}, it was proven that the mass loss due to the penalty term is of order~$\mathcal{O}(\varepsilon^{\frac12})$; however, 
an order~$\mathcal{O}(\varepsilon)$ was numerically observed.
Therefore, as the introduction of the penalty 
term increases the error of the method,
it should be used only when it is necessary to get convergence of the non-linear Newton solver (see Remark~\ref{RMK:PENALTY-TERM}). 
For this reason, in the next steps of this test case, we 
fix~$\varepsilon = 0$, 
as  the 
penalty term does not seem to be necessary. 
\subsubsection*{Convergence with respect to the polynomial degree~$\ell$}
\begin{figure}[t]
	\centering
	{\includegraphics[width=\textwidth]{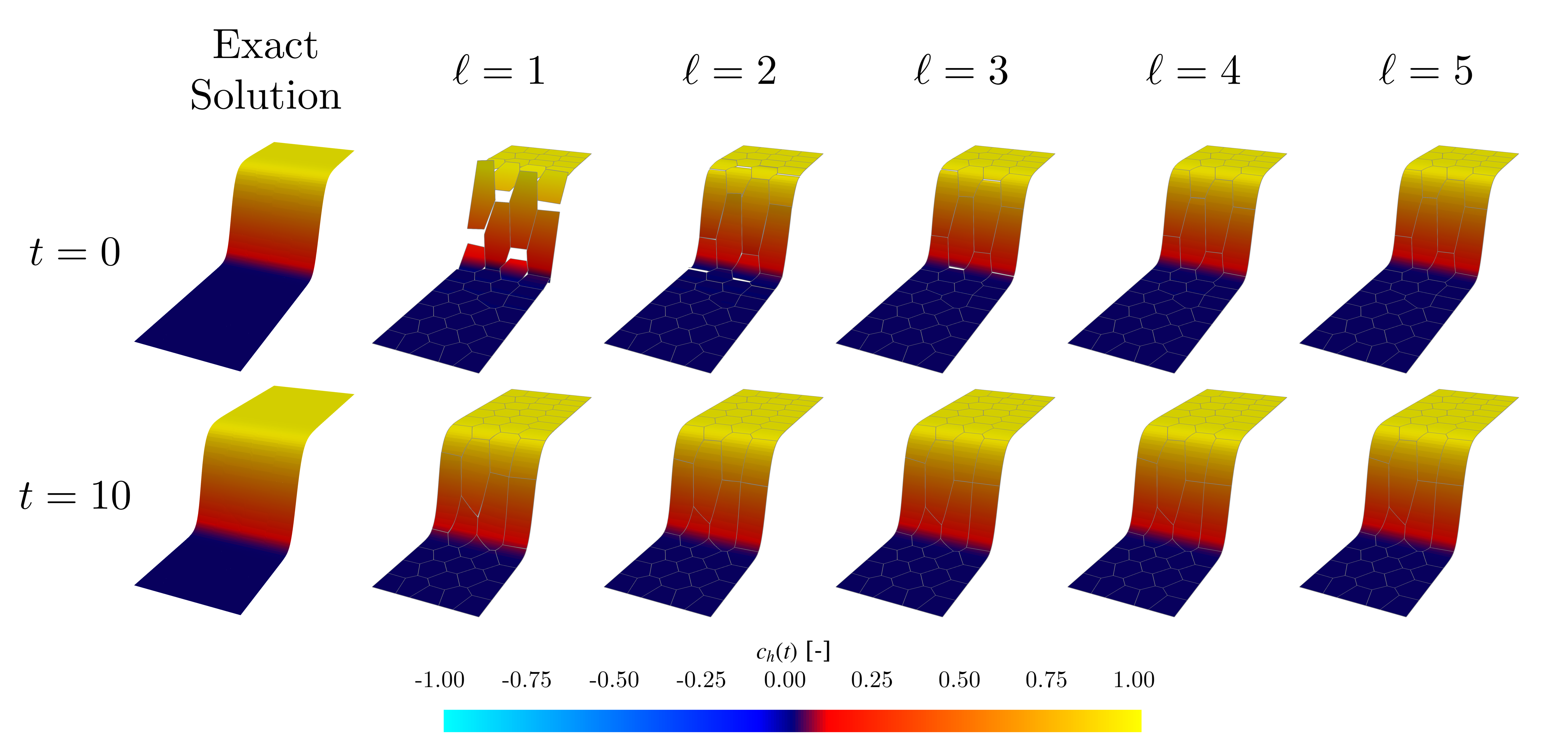}}
	\caption{Test case 2: Comparison of initial conditions and solutions of wave propagation problem for different polynomial degrees $\ell=1,...,5$ at different times $t=0,10$.}
	\label{fig:waves2D}
\end{figure}
In the second test, we analyze the convergence with respect to the polynomial degree~$\ell$ for different choices of \BDF{\bdfp} scheme~($\bdfp=1,...,6$). We use a fixed mesh consisting of 50 elements generated with PolyMesher \cite{talischi_polymesher_2012}, 
and consider as time step $\tau = 2.5\times 10^{-2}$ and as final time $T=10$.
\par
We report in Figure \ref{fig:waves2D} the numerical solutions obtained with different choices of the polynomial degree~$\ell$, $\ell=1,...,5$, and with the \BDF{6} scheme for the time discretization.  On the first row, we report the initial condition at time $t=0$. 
As shown, the choice of using an~$L^2(\Omega)$-projection of the field~$c_0(\boldsymbol{x})$ can 
result in a highly discontinuous 
approximation and. Locally, this can cause the values to fall outside the range~$(0,1)$, particularly for low polynomial degrees. However, the transformation in~\eqref{EQ::UPUQ} ensures that the solution values are pointwise within~$(0,1)$ from the very first time step. 
As shown in the second row, the method provides an accurate representation of the wavefront propagation.
\par
In Figure \ref{fig:errorswaves_ell}, we report the errors~$\Ec$ and~$\Esigma$.
We observe that using a higher-order time-stepping scheme results in spectral convergence in~$\ell$. In contrast, the \BDF{1} and \BDF{2} schemes 
lead to error stagnation, which cannot be reduced by improving the space discretization, as the time discretization error dominates in such cases.
We do not report the errors of \BDF{\bdfp} methods for \bdfp=$4,5,6$ as their performance in this test case is identical to that of the \BDF{3} scheme.

\subsubsection*{Convergence with respect to the space mesh size $h$}
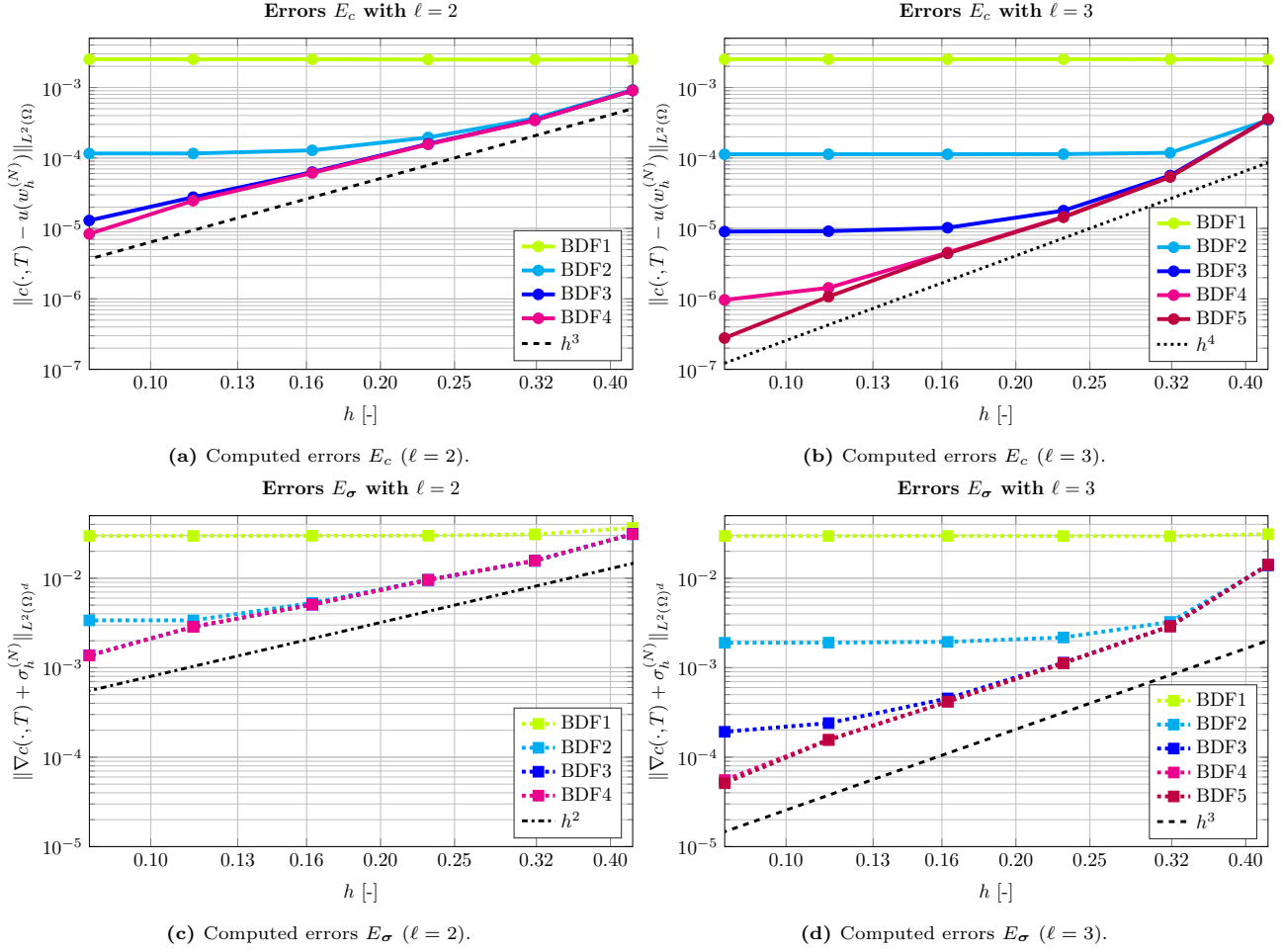
\begin{figure}[t!]
    \begin{subfigure}[b]{0.5\textwidth}
          \resizebox{\textwidth}{!}{\begin{tikzpicture}

\begin{axis}[%
width=3.875in,
height=2.36in,
at={(2.6in,1.099in)},
scale only axis,
xmode=log,
xmin=0.0831,
xmax=0.4278,
xminorticks=true,
xlabel = {$h$ [-]},
ylabel = {$\|c(\cdot,T)-u(w_h^{(N)})\|_{L^2(\Omega)}$},
xtick={0.10, 0.13, 0.16, 0.20, 0.25, 0.32, 0.4},
xticklabel={\pgfmathparse{exp(\tick)}\pgfmathprintnumber{\pgfmathresult}},
x tick label style={
/pgf/number format/.cd, fixed, fixed zerofill,
precision=2},
ymode=log,
ymin=1e-7,
ymax=5e-3,
yminorticks=true,
axis background/.style={fill=white},
title style={font=\bfseries},
title={Errors $E_c$ with $\ell=2$},
xmajorgrids,
xminorgrids,
ymajorgrids,
yminorgrids,
legend style={at={(0.98,0.42)},legend cell align=left, draw=white!15!black}
]
              
\addplot [color=lime, line width=2.0pt, mark=*]
  table[row sep=crcr]{%
0.4278  2.5000e-03 \\
0.3185  2.4902e-03 \\
0.2309  2.4966e-03 \\
0.1628  2.5162e-03 \\
0.1137 	2.5189e-03 \\
0.0831  2.5210e-03 \\
};
\addlegendentry{\BDF{1}}

\addplot [color=cyan, line width=2.0pt, mark=*]
  table[row sep=crcr]{%
0.4278  9.2478e-04 \\
0.3185  3.6483e-04 \\
0.2309  1.9567e-04 \\
0.1628  1.2843e-04 \\
0.1137 	1.1578e-04 \\
0.0831  1.1578e-04 \\
};
\addlegendentry{\BDF{2}}

\addplot [color=blue, line width=2.0pt, mark=*]
  table[row sep=crcr]{%
0.4278  9.0899e-04 \\
0.3185  3.4420e-04 \\
0.2309  1.5844e-04 \\
0.1628  6.2963e-05 \\
0.1137 	2.7753e-05 \\
0.0831  1.2949e-05 \\
};
\addlegendentry{\BDF{3}}

\addplot [color=magenta, line width=2.0pt, mark=*]
  table[row sep=crcr]{%
0.4278  9.0725e-04 \\
0.3185  3.4296e-04 \\
0.2309  1.5612e-04 \\
0.1628  6.1298e-05 \\
0.1137 	2.4744e-05 \\
0.0831  8.3508e-06 \\
};
\addlegendentry{\BDF{4}}

\addplot [color=black, line width=1.5pt, dashed]
  table[row sep=crcr]{%
0.5000  8.0000e-04 \\
0.1000  6.4000e-06 \\
0.0200  5.1200e-08 \\
};
\addlegendentry{$h^3$}

\end{axis}
\end{tikzpicture}%
}
        \caption{Computed errors 
        $\Ec$ $(\ell=2)$.}
        \label{fig:errorswavesh_l2_L2}
    \end{subfigure}%
    \begin{subfigure}[b]{0.5\textwidth}
        \resizebox{\textwidth}{!}{\begin{tikzpicture}

\begin{axis}[%
width=3.875in,
height=2.36in,
at={(2.6in,1.099in)},
scale only axis,
xmode=log,
xmin=0.0831,
xmax=0.4278,
xminorticks=true,
xlabel = {$h$ [-]},
ylabel = {$\|c(\cdot,T)-u(w_h^{(N)})\|_{L^2(\Omega)}$},
xtick={0.10, 0.13, 0.16, 0.20, 0.25, 0.32, 0.4},
xticklabel={\pgfmathparse{exp(\tick)}\pgfmathprintnumber{\pgfmathresult}},
x tick label style={
/pgf/number format/.cd, fixed, fixed zerofill, precision=2},
ymode=log,
ymin=1e-7,
ymax=5e-3,
yminorticks=true,
axis background/.style={fill=white},
title style={font=\bfseries},
title={Errors $E_c$ with $\ell=3$},
xmajorgrids,
xminorgrids,
ymajorgrids,
yminorgrids,
legend style={at={(0.98,0.49)},legend cell align=left, draw=white!15!black}
]
              
\addplot [color=lime, line width=2.0pt, mark=*]
  table[row sep=crcr]{%
0.4278  2.5000e-03 \\
0.3185  2.5139e-03 \\
0.2309  2.5206e-03 \\
0.1628  2.5215e-03 \\
0.1137 	2.5215e-03 \\
0.0831  2.5216e-03 \\
};
\addlegendentry{\BDF{1}}

\addplot [color=cyan, line width=2.0pt, mark=*]
  table[row sep=crcr]{%
0.4278  3.4487e-04 \\
0.3185  1.1846e-04 \\
0.2309  1.1338e-04 \\
0.1628  1.1295e-04 \\
0.1137 	1.1287e-04 \\
0.0831  1.1286e-04 \\
};
\addlegendentry{\BDF{2}}

\addplot [color=blue, line width=2.0pt, mark=*]
  table[row sep=crcr]{%
0.4278  3.5558e-04 \\
0.3185  5.6385e-05 \\
0.2309  1.7888e-05 \\
0.1628  1.0227e-05 \\
0.1137 	9.1235e-06 \\
0.0831  9.0256e-06 \\
};
\addlegendentry{\BDF{3}}

\addplot [color=magenta, line width=2.0pt, mark=*]
  table[row sep=crcr]{%
0.4278  3.5719e-04 \\
0.3185  5.4155e-05 \\
0.2309  1.4585e-05 \\
0.1628  4.5461e-06 \\
0.1137 	1.4352e-06 \\
0.0831  9.6626e-07 \\
};
\addlegendentry{\BDF{4}}

\addplot [color=purple, line width=2.0pt, mark=*]
  table[row sep=crcr]{%
0.4278  3.5739e-04 \\
0.3185  5.3812e-05 \\
0.2309  1.4411e-05 \\
0.1628  4.4206e-06 \\
0.1137 	1.0751e-06 \\
0.0831  2.7815e-07 \\
};
\addlegendentry{\BDF{5}}
					
\addplot [color=black, line width=1.5pt, dotted]
  table[row sep=crcr]{%
0.5000  1.6000e-04 \\
0.1000  2.5600e-07 \\
0.0200  4.0960e-10 \\
};
\addlegendentry{$h^4$}

\end{axis}
\end{tikzpicture}%
}
        \caption{Computed errors 
        $\Ec$ $(\ell=3)$.}
        \label{fig:errorswavesh_l3_L2}
    \end{subfigure}\\[0.5mm]
    \begin{subfigure}[b]{0.5\textwidth}
          \resizebox{\textwidth}{!}{\begin{tikzpicture}

\begin{axis}[%
width=3.875in,
height=2.36in,
at={(2.6in,1.099in)},
scale only axis,
xmode=log,
xmin=0.0831,
xmax=0.4278,
xminorticks=true,
xlabel = {$h$ [-]},
ylabel = {$\|\nabla c(\cdot,T) + \sigma_h^{(N)}\|_{L^2(\Omega)^d}$},
xtick={0.10, 0.13, 0.16, 0.20, 0.25, 0.32, 0.4},
xticklabel={\pgfmathparse{exp(\tick)}\pgfmathprintnumber{\pgfmathresult}},
x tick label style={
/pgf/number format/.cd, fixed, fixed zerofill,
precision=2},
ymode=log,
ymin=1e-5,
ymax=5e-2,
yminorticks=true,
axis background/.style={fill=white},
title style={font=\bfseries},
title={Errors $E_{\boldsymbol{\sigma}}$ with $\ell=2$},
xmajorgrids,
xminorgrids,
ymajorgrids,
yminorgrids,
legend style={at={(0.98,0.42)},legend cell align=left, draw=white!15!black}
]
              
\addplot [color=lime, line width=2.0pt, mark=square*, dotted, mark options=solid]
  table[row sep=crcr]{%
0.4278  3.6400e-02 \\					
0.3185  3.0918e-02 \\
0.2309  2.9923e-02 \\
0.1628  2.9845e-02 \\
0.1137 	2.9687e-02 \\
0.0831  2.9649e-02 \\
};
\addlegendentry{\BDF{1}}

\addplot [color=cyan, line width=2.0pt, mark=square*, dotted, mark options=solid]
  table[row sep=crcr]{%
0.4278  3.1166e-02 \\ 	 	 			
0.3185  1.5600e-02 \\
0.2309  9.5889e-03 \\
0.1628  5.3003e-03 \\
0.1137 	3.3811e-03 \\
0.0831  3.3811e-03 \\
};
\addlegendentry{\BDF{2}}

\addplot [color=blue, line width=2.0pt, mark=square*, dotted, mark options=solid]
  table[row sep=crcr]{%
0.4278  3.1146e-02 \\ 	 				
0.3185  1.5680e-02 \\
0.2309  9.5778e-03 \\
0.1628  5.0701e-03 \\
0.1137 	2.8646e-03 \\
0.0831  1.3758e-03 \\
};
\addlegendentry{\BDF{3}}

\addplot [color=magenta, line width=2.0pt, mark=square*, dotted, mark options=solid]
  table[row sep=crcr]{%
0.4278  3.1146e-02 \\
0.3185  1.5690e-02 \\
0.2309  9.5798e-03 \\
0.1628  5.0767e-03 \\
0.1137 	2.8618e-03 \\
0.0831  1.3675e-03 \\
};
\addlegendentry{\BDF{4}}

\addplot [color=black, line width=1.5pt, dashdotted]
  table[row sep=crcr]{%
0.5000  2.0000e-02 \\
0.1000  8.0000e-04 \\
0.0200  3.2000e-05 \\
};
\addlegendentry{$h^2$}

\end{axis}
\end{tikzpicture}%
}
        \caption{Computed errors 
        $\Esigma$ $(\ell=2)$.}
        \label{fig:errorswavesh_l2_dG}
    \end{subfigure}%
    \begin{subfigure}[b]{0.5\textwidth}
        \resizebox{\textwidth}{!}{\begin{tikzpicture}

\begin{axis}[%
width=3.875in,
height=2.36in,
at={(2.6in,1.099in)},
scale only axis,
xmode=log,
xmin=0.0831,
xmax=0.4278,
xminorticks=true,
xlabel = {$h$ [-]},
ylabel = {$\|\nabla c(\cdot,T) + \sigma_h^{(N)}\|_{L^2(\Omega)^d}$},
xtick={0.10, 0.13, 0.16, 0.20, 0.25, 0.32, 0.4},
xticklabel={\pgfmathparse{exp(\tick)}\pgfmathprintnumber{\pgfmathresult}},
x tick label style={
/pgf/number format/.cd, fixed, fixed zerofill,
precision=2},
ymode=log,
ymin=1e-5,
ymax=5e-2,
yminorticks=true,
axis background/.style={fill=white},
title style={font=\bfseries},
title={Errors $E_{\boldsymbol{\sigma}}$ with $\ell=3$},
xmajorgrids,
xminorgrids,
ymajorgrids,
yminorgrids,
legend style={at={(0.98,0.49)},legend cell align=left, draw=white!15!black}
]
              
\addplot [color=lime, line width=2.0pt, mark=square*, dotted, mark options=solid]
  table[row sep=crcr]{%
0.4278  3.1100e-02 \\
0.3185  2.9432e-02 \\
0.2309  2.9608e-02 \\
0.1628  2.9610e-02 \\
0.1137 	2.9605e-02 \\
0.0831  2.9603e-02 \\
};
\addlegendentry{\BDF{1}}

\addplot [color=cyan, line width=2.0pt, mark=square*, dotted, mark options=solid]
  table[row sep=crcr]{%
0.4278  1.3800e-02 \\
0.3185  3.2445e-03 \\
0.2309  2.1753e-03 \\
0.1628  1.9438e-03 \\
0.1137 	1.9054e-03 \\
0.0831  1.8985e-03 \\
};
\addlegendentry{\BDF{2}}

\addplot [color=blue, line width=2.0pt, mark=square*, dotted, mark options=solid]
  table[row sep=crcr]{%
0.4278  1.4100e-02 \\ 					
0.3185  2.9110e-03 \\
0.2309  1.1370e-03 \\
0.1628  4.5338e-04 \\
0.1137 	2.4009e-04 \\
0.0831  1.9211e-04 \\
};
\addlegendentry{\BDF{3}}

\addplot [color=magenta, line width=2.0pt, mark=square*, dotted, mark options=solid]
  table[row sep=crcr]{%
0.4278  1.4100e-02 \\					
0.3185  2.8999e-03 \\
0.2309  1.1253e-03 \\
0.1628  4.1746e-04 \\
0.1137 	1.5663e-04 \\
0.0831  5.5498e-05 \\
};
\addlegendentry{\BDF{4}}

\addplot [color=purple, line width=2.0pt, mark=square*, dotted, mark options=solid]
  table[row sep=crcr]{%
0.4278  1.4100e-02 \\					
0.3185  2.8953e-03 \\				
0.2309  1.1231e-03 \\
0.1628  4.1704e-04 \\
0.1137 	1.5544e-04 \\
0.0831  5.1187e-05 \\
};
\addlegendentry{\BDF{5}}
					
\addplot [color=black, line width=1.5pt, dashed]
  table[row sep=crcr]{%
0.5000  3.2000e-03 \\
0.1000  2.5600e-05 \\
0.0200  2.0480e-07 \\
};
\addlegendentry{$h^3$}

\end{axis}
\end{tikzpicture}%
}
        \caption{Computed errors 
        $\Esigma$ $(\ell=3)$.}
        \label{fig:errorswavesh_l3_dG}
    \end{subfigure}%
    \caption{Test case 2: computed errors and convergence rates.
    }
    \label{fig:errorswavesh}
\end{figure}
In the third test, we analyze the convergence with respect to the space mesh size~$h$ for different choices of the \BDF{\bdfp} scheme ($\bdfp=1,...,6$). 
We fix two different polynomial degrees of approximation,~$\ell=2$ and~$\ell=3$, and test on a sequence of meshes consisting of~$50, 100, 200, 400, 800, 1600$ elements generated with PolyMesher~\cite{talischi_polymesher_2012}. 
We take as final time~$T = 1$, and a fixed time step~$\tau = 10^{-1}$ and use high-order~\BDF{\bdfp} time stepping schemes to let the space error dominate.
\par
In Figures \ref{fig:errorswavesh_l2_L2} and \ref{fig:errorswavesh_l2_dG}, we report the errors~$\Ec$ and~$\Esigma$, respectively,
for~$\ell=2$. Additionally, in Figures \ref{fig:errorswavesh_l3_L2} and \ref{fig:errorswavesh_l3_dG}, we report the errors~$\Ec$ and~$\Esigma$ 
for~$\ell=3$. For sufficiently high-order time stepping schemes,
we can observe that the 
the errors~$\Ec$ and~$\Esigma$ decay with optimal orders~$\mathcal{O}(h^{\ell+1})$ and~$\mathcal{O}(h^{\ell})$, respectively. 
In contrast, low-order BDF schemes cause error stagnation, as the time discretization error dominates over the error in space.
For the chosen mesh refinements, the correct slope can be found using the \BDF{4} scheme for $\ell=2$, but it requires a \BDF{5} for $\ell=3$. In Figure \ref{fig:errorswavesh}, we do not report the results for \BDF{\bdfp} schemes of higher order, as they do not improve significantly the accuracy of the method. This experiment also emphasizes the \emph{unconditional stability} of the method, meaning that there is no need to choose a time step sufficiently small with respect to the space mesh size to prevent blow-up of the numerical solution.

\subsubsection*{Comparison with an interior penalty DG method}
Finally, we compare the results obtained with our method and the ones obtained with another method proposed in the literature. In particular, we focus on an IPDG method proposed in~\cite{Corti_Bonizzoni_Dede_Quarteroni_Antonietti:2023}.
This method is not positivity-preserving, but if the mesh size is sufficiently fine or the polynomial degree is sufficiently high, the analytical solution is approximated correctly.
\par
For this test, we use two meshes with 50 and 200 elements and different polynomial degrees~$\ell=1,...,5$. We consider as time step~$\tau = 2.5\times 10^{-2}$ and as final time~$T=10$. 
As for the time discretization, we adopt the \BDF{1} and \BDF{2} schemes for our structure-preserving LDG method. The results are compared with those obtained with the implicit Euler 
scheme (\BDF{1}) and the Crank-Nicolson scheme (of order 2) for the IPDG approach proposed in~\cite{Corti_Bonizzoni_Dede_Quarteroni_Antonietti:2023}. This ensures that, in the comparison, the two space discretizations are paired with time stepping schemes of the same orders.
\par
\begin{figure}[t]
	\centering
	{\includegraphics[width=\textwidth]{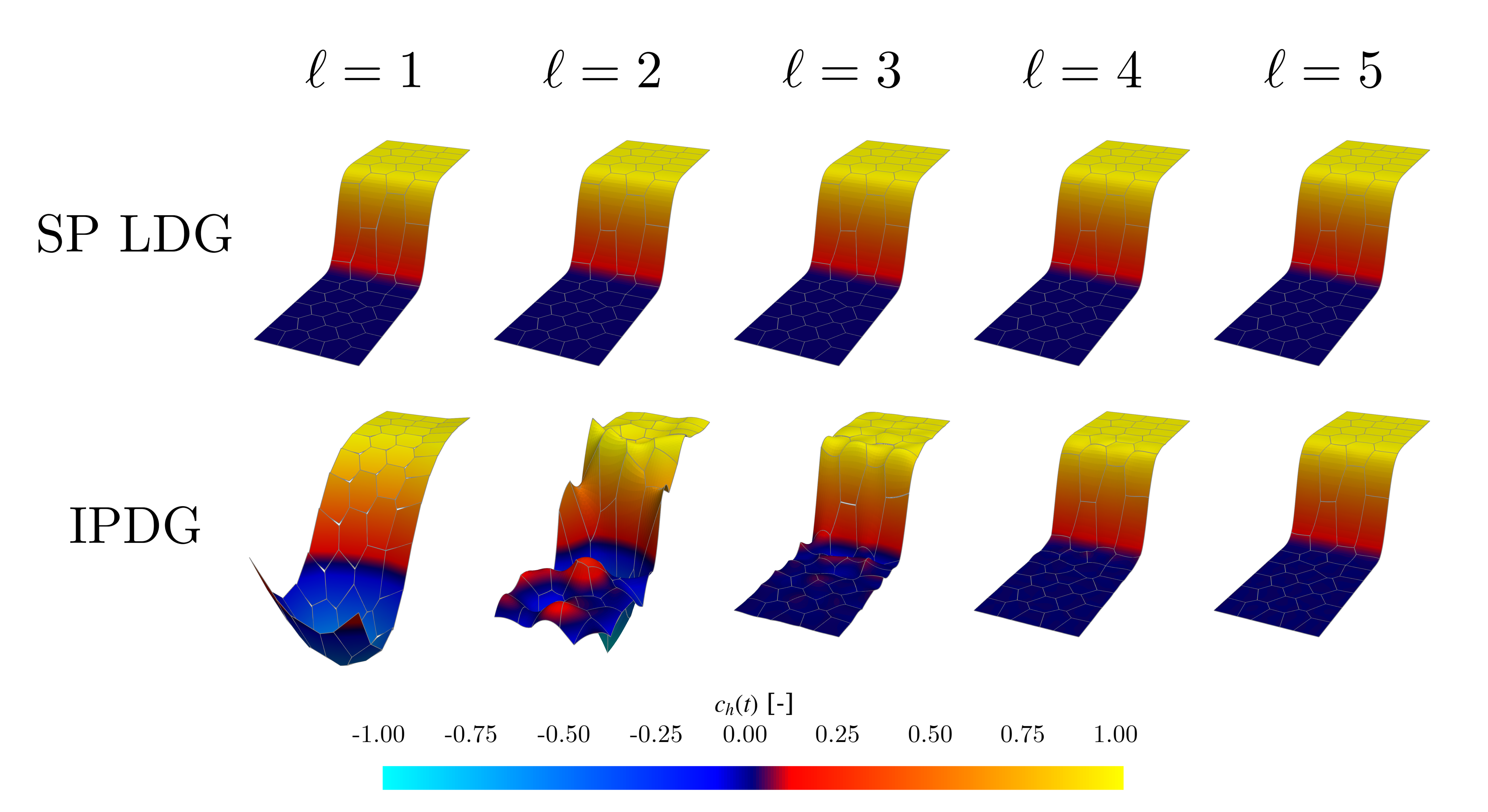}}
	\caption{Test case 2: Comparison solutions at time $t=3$ for different polynomial degrees $\ell=1,...,5$ and computed with different methods: the structure-preserving LDG (first line) and the IPDG (second line).}
	\label{fig:L-IP-DGcomparison}
\end{figure}
In Figure \ref{fig:L-IP-DGcomparison}, we report the results obtained for a mesh with 50 elements, using the time integration schemes of the second-order. 
We observe that, using the IPDG method, the solution does not maintain positivity and takes negative values at certain points, thereby losing physical meaning. This issue can be mitigated by increasing the degree~$\ell$ of the approximation in space. In contrast, our method remains free of oscillations and preserves the pointwise bounds of the continuous solution.
\par
\begin{table}[t]
    \centering
    \setlength{\extrarowheight}{1.5pt}
    \begin{tabular}{|c|c|c c c c c|}
    \multicolumn{7}{c}{$\boldsymbol{h \approx 0.4278}$ \textbf{and} $\boldsymbol{\tau = 2.5\times10^{-2}}$} \\[1pt]
    \hline
    \multicolumn{2}{|c|}{\textbf{Method}} & $\ell=1$ & $\ell=2$ & $\ell=3$ & $\ell=4$ & $\ell=5$
    \\ \hline 
    \multirow{2}{*}{\BDF{1}} & SP-LDG 
    & $4.72\times10^{-2}$ & $2.28\times10^{-2}$  
    & $1.50\times10^{-2}$ & $1.29\times10^{-2}$ 
    & $1.43\times10^{-2}$ 
    \\ 
    & IPDG 
    & \color{red}$\mathbf{1.06}$ & \color{red}$\mathbf{5.98}$  
    & \color{red}$\mathbf{3.33\times10^{2}}$ & $1.43\times10^{-2}$ 
    & $1.43\times10^{-2}$ 
    \\ \hline   
    \BDF{2} & SP-LDG 
    & $4.72\times10^{-2}$ & $3.73\times10^{-3}$  
    & $1.13\times10^{-3}$ & $2.84\times10^{-4}$ 
    & $2.50\times10^{-4}$ 
    \\ CN
    & IPDG 
    & \color{red}$\mathbf{3.38\times10^{3}}$ & \color{red}$\mathbf{1.90\times10^{3}}$  
    & \color{red}$\mathbf{2.33\times10^{1}}$ & $1.15\times10^{-1}$ 
    & $7.68\times10^{-3}$ 
    \\ \hline
    \end{tabular}
    \\[4pt]
    \begin{tabular}{|c|c|c c c c c|}
    \multicolumn{7}{c}{$\boldsymbol{h \approx 0.2309}$ \textbf{and} $\boldsymbol{\tau = 2.5\times10^{-2}}$} \\[1pt]
    \hline
    \multicolumn{2}{|c|}{\textbf{Method}} & $\ell=1$ & $\ell=2$ & $\ell=3$ & $\ell=4$ & $\ell=5$
    \\ \hline 
    \multirow{2}{*}{\BDF{1}} & SP-LDG 
    & $4.67\times10^{-2}$ & $1.31\times10^{-2}$  
    & $1.43\times10^{-2}$ & $1.43\times10^{-2}$ 
    & $1.02\times10^{-2}$ 
    \\ 
    & IPDG 
    & \color{red}$\mathbf{1.05}$ & \color{red}$\mathbf{1.07}$  
    & $4.32\times10^{-2}$ & $1.66\times10^{-2}$ 
    & $1.69\times10^{-2}$ 
    \\ \hline   
    \BDF{2} & SP-LDG 
    & $3.18\times10^{-2}$ & $1.44\times10^{-3}$  
    & $2.46\times10^{-4}$ & $2.52\times10^{-4}$ 
    & $2.84\times10^{-4}$ 
    \\ CN
    & IPDG 
    & \color{red}$\mathbf{3.13\times10^{3}}$  & \color{red}$\mathbf{9.66\times10^{4}}$  
    & $3.55\times10^{-2}$ & $6.17\times10^{-4}$ 
    & $9.21\times10^{-5}$ 
    \\ \hline
    \end{tabular}
    \\[4pt]
    \begin{tabular}{|c|c|c c c c c|}
    \multicolumn{7}{c}{$\boldsymbol{h \approx 0.4278}$ \textbf{and} $\boldsymbol{\tau = 5.0\times10^{-3}}$} \\[1pt]
    \hline
    \multicolumn{2}{|c|}{\textbf{Method}} & $\ell=1$ & $\ell=2$ & $\ell=3$ & $\ell=4$ & $\ell=5$
    \\ \hline 
    \multirow{2}{*}{\BDF{1}} & SP-LDG 
    & $3.66\times10^{-2}$ & $5.34\times10^{-3}$  
    & $4.91\times10^{-3}$ & $4.20\times10^{-3}$ 
    & $4.25\times10^{-3}$ 
    \\
    & IPDG 
    & \color{red}$\mathbf{1.37}$ & \color{red}$\mathbf{2.89}$  
    & \color{red}$\mathbf{2.62\times10^{3}}$ & $4.51\times10^{-3}$ 
    & $2.97\times10^{-3}$ 
    \\ \hline 
    \BDF{2} & SP-LDG 
    & $3.14\times10^{-2}$ & $4.48\times10^{-3}$  
    & $8.59\times10^{-4}$ & $1.15\times10^{-4}$ 
    & $2.55\times10^{-5}$ 
    \\ CN
    & IPDG 
    & \color{red}$\mathbf{2.51\times10^{4}}$ & \color{red}$\mathbf{2.60\times10^{4}}$  
    & \color{red}$\mathbf{1.59\times10^{4}}$ & $4.47\times10^{-4}$ 
    & $2.77\times10^{-4}$ 
    \\ \hline
    \end{tabular}
    \\[4pt]
    \begin{tabular}{|c|c|c c c c c|}
    \multicolumn{7}{c}{$\boldsymbol{h \approx 0.2309}$ \textbf{and} $\boldsymbol{\tau = 5.0\times10^{-3}}$} \\[1pt]
    \hline
    \multicolumn{2}{|c|}{\textbf{Method}} & $\ell=1$ & $\ell=2$ & $\ell=3$ & $\ell=4$ & $\ell=5$
    \\ \hline 
    \multirow{2}{*}{\BDF{1}} & SP-LDG 
    & $3.46\times10^{-2}$ & $1.97\times10^{-3}$  
    & $2.81\times10^{-3}$ & $2.82\times10^{-3}$ 
    & $2.83\times10^{-3}$
    \\
    & IPDG 
    & \color{red}$\mathbf{1.08}$ & \color{red}$\mathbf{1.09}$  
    & $3.65\times10^{-2}$ & $3.05\times10^{-3}$ 
    & $3.33\times10^{-3}$ 
    \\ \hline 
    \BDF{2} & SP-LDG 
    & $3.17\times10^{-2}$ & $1.62\times10^{-3}$  
    & $1.54\times10^{-5}$ & $1.05\times10^{-5}$ 
    & $1.03\times10^{-5}$
    \\ CN
    & IPDG 
    & \color{red}$\mathbf{2.75\times10^{4}}$ & \color{red}$\mathbf{3.39\times10^{4}}$  
    & $3.55\times10^{-2}$ & $6.57\times10^{-4}$ 
    & $3.97\times10^{-5}$
    \\ \hline
    \end{tabular}
    \caption{Computed errors in the $L^2(\Omega)$ norm at final time $T=10$ with different methods (SP-LDG and IPDG~\cite{Corti_Bonizzoni_Dede_Quarteroni_Antonietti:2023}).
    }    \label{tab:errorspconvergence}
\end{table}
To quantify the quality of the numerical approximation obtained with the two methods, we report in Table~\ref{tab:errorspconvergence} the errors in the~$L^2(\Omega)$ norm, namely~$\Norm{c(\cdot,T)-u(w_h^{(N)})}{L^2(\Omega)}$ for the structure-preserving LDG method and~$\Norm{c(\cdot,T)-c_h^{(N)}}{L^2(\Omega)}$ for the IPDG method. 
This shows that the method proposed in \cite{Corti_Bonizzoni_Dede_Quarteroni_Antonietti:2023} only approximates the exact solution accurately when a sufficiently large number of degrees of freedom are used. This requires increasing the accuracy in space (by either increasing~$\ell$ or reducing~$h$), which appears to be more crucial than reducing the time step. 
In contrast, our method is not affected by these issues and can capture the solution correctly even with low-order polynomial degrees. The performance of both methods becomes comparable for higher-order approximations.
\section{Spreading of \textalpha{}-synuclein in a two-dimensional brain section}\label{sec:brain}

In this section, we present some numerical simulations of the spreading of \textalpha{}-synuclein, conducted using our structure-preserving LDG method. The aggregation, phosphorylation, and nitration of this synaptic protein have been suggested to be critical processes in forming aggregates known as Lewy bodies. These lesions are associated with several neurodegenerative diseases, including Parkinson's disease (PD), dementia with Lewy bodies (DLB), incidental Lewy body disease (ILBD), and Alzheimer’s disease with Lewy bodies (ADLB) \cite{beach_multi-organ_2010}. All these conditions share a common developmental process described in the literature by the Unified Staging System for Lewy Body Disorders (USSLB) \cite{beach_unified_2009,adler_unified_2019}. This theory divides pathological spreading of \textalpha{}-synuclein into 4 stages, but with two distinct patterns that may emerge at the prodromal stage: primarily brainstem inclusions (most common in PD) or limbic inclusions (typical in DLB, ADLB)  \cite{adler_unified_2019}.
\par
In the literature, the mathematical modeling of these processes typically uses the \emph{heterodimer model} ~\cite{fornari_prion-like_2019,weickenmeierPhysicsbasedModelExplains2019,thompson_proteinprotein_2020,Antonietti_Bonizzoni_Corti_Dallolio:2024}, which describes the evolution of the healthy protein concentration~$p$ and the misfolded protein concentration~$q$. Denoting as~$\kappa_p>0$ the production rate of~$p$, $\lambda_p>0$ and~$\lambda_q> 0$ the clearance rates of~$p$ and~$q$, respectively, and~$\mu_{pq}>0$ the conversion rate from~$p$ to~$q$, the heterodimer system with initial conditions and open boundary conditions reads as follows:
\begin{alignat*}{3}
\dpt p - \nabla \cdot (\D \nabla p) & = - p(\lambda_p + \mu_{pq} q) + \kappa_p & & \quad \text{ in } \QT,\\
\dpt q - \nabla \cdot (\D \nabla q) & = -q(\lambda_q - \mu_{pq} p) & & \quad \text{ in } \QT,\\
(\D \nabla p) \cdot \bnOmega = 0  \ \text{ and } \ (\D \nabla q) \cdot \bnOmega & = 0 & & \quad  \text{ on } \Gamma \times (0, T),\\
p(\cdot, 0) = p_{0} \text{ and } 
q(\cdot, 0) & = q_{0} & & \quad \text{ in } \Omega.
\end{alignat*}
Assuming constant coefficients and ~$p\gg q$, and neglecting the time derivative and diffusion of~$p$, $p$ can be computed in terms of~$q$ as~$p=\frac{\kappa_p}{\lambda_p+\mu_{pq} q}$. By performing a Taylor expansion,~$p\simeq \frac{\kappa_p}{\lambda_p}\left(1-\frac{\mu_{pq}}{\lambda_p}q\right)$. Substituting this into the second equation of the heterodimer system, and setting~$c:=q/q_M$, with~$q_M:=\lambda_p(\kappa_p \mu_{pq} - \lambda_p\lambda_q)/(\kappa_p\mu_{pq}^2)$, one obtains the FK equation with~$\alpha=(\kappa_p \mu_{pq} - \lambda_p\lambda_q)/{\lambda_p}$; see~\cite{fornari_prion-like_2019}. Under these assumptions, the FK equation in \eqref{EQN::FISHER-KPP} can be viewed as a simplified version of the heterodimer model, applicable when the dynamics of the healthy protein population are negligible. Indeed, the FK equation has been widely used in the literature to model and simulate the spreading of proteins in neurodegenerative diseases \cite{weickenmeierPhysicsbasedModelExplains2019, Corti_Bonizzoni_Dede_Quarteroni_Antonietti:2023}.
\par
\begin{figure}[t!]
    \centering
    {\includegraphics[width=\textwidth]{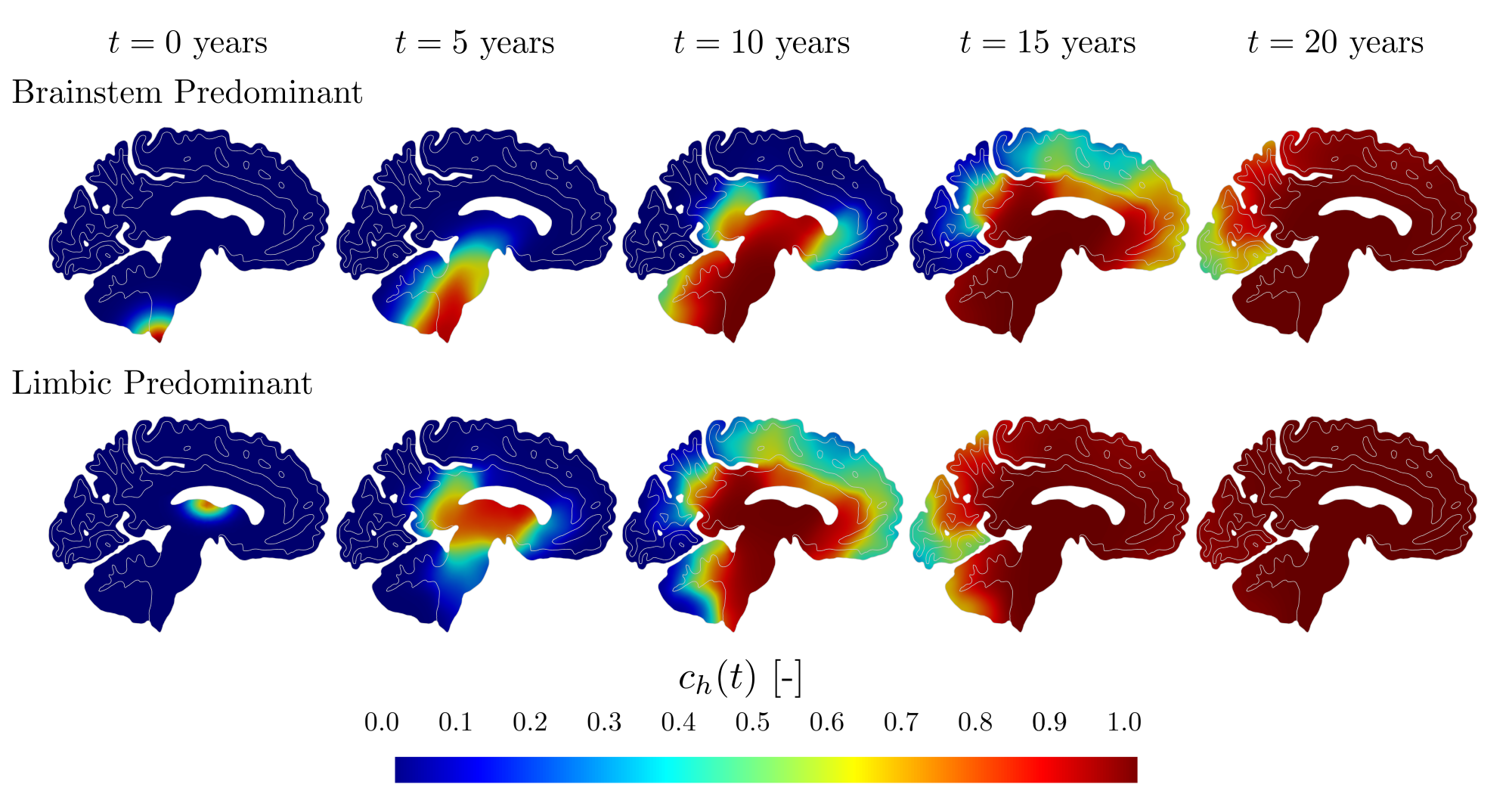}}
    \caption{Test case 3: Patterns of \textalpha{}-synuclein concentration at different stages of the pathology.} 
    \label{fig:Solution2DPark}
\end{figure}

Here, we simulate the spreading of \textalpha{}-synuclein in both prodromal stages, modeled by the FK equation to analyze the evolution in a simplified 2D geometry. The brain geometry is segmented starting from a structural magnetic resonance image of a brain from the OASIS-3 database \cite{lamontagne_oasis-3_2019}. We construct a polytopal agglomerated grid of a sagittal 2D brain section. The initial mesh of the brain slice is a triangular mesh composed of $43\,402$ elements constructed using VMTK \cite{antiga_image-based_2008}. As presented in \cite[\S5.1]{Corti_Bonizzoni_Antonietti:2024}, the final mesh is a polytopal mesh consisting of $534$ elements constructed by agglomerating the elements of the initial triangular mesh. The use of a DG method allows for meshes with elements of fairly general shape, while the construction of agglomerated elements helps to reduce the computational cost. Each element is labeled as either white or grey matter, according to the segmentation of magnetic resonance images. Another advantage of mesh agglomeration is that it allows the external boundary and the internal interface between white and grey matter to be described using just a few polytopal elements.
\par
In the simulations, we set
the physical parameters as follows~\cite{schafer_interplay_2019,weickenmeierPhysicsbasedModelExplains2019}: 
the reaction coefficient is
\[
\alpha=\alpha(\bx)=\begin{cases}
 0.45\;\mathrm{year}^{-1} & \text{ in the grey matter}, \\
 0.9\;\mathrm{year}^{-1} & \text{ in the white matter},
\end{cases}
\]
and the diffusion tensor, which is isotropic in the grey matter and anisotropic in the white matter, is
\[
\D=\D(\bx)=\begin{cases}
 \dext\mathbb{I}_2 & \text{ in the grey matter}, \\
 \dext\mathbb{I}_2 + \daxn\, \ba(\bx) \otimes \ba(\bx) & \text{ in the white matter},
\end{cases}
\]
with isotropic diffusion with the same coefficient~$\dext = 8\;\mathrm{mm}^2/\mathrm{year}$, and additional anisotropic diffusion in the 
white matter 10 times faster, namely~$\daxn = 80\;\mathrm{mm}^2/\mathrm{year}$. The latter is associated with axonal directions $\ba=\ba(\bx)$, with~$|\ba|=1$, derived from diffusion-weighted images as in \cite{Corti_Bonizzoni_Dede_Quarteroni_Antonietti:2023}.
The diffusion tensor~$\D$ is uniformly positive definite in the domain, with constant larger than or equal to~$\dext$.
\par
Concerning the numerical discretization of the problem, we use  
polynomials of degree~$\ell = 2$, the power-mean parameter~$\theta = 1/2$ in~\eqref{EQN::DEF-h}, and a stabilization parameter~$\eta_0 = 2$ for the space discretization. For the time discretization, we choose a time step $\tau = 2.5\times 10^{-2}\,\mathrm{year}$, a final time $T=25\;\mathrm{years}$, and we use the \BDF{6} time-stepping scheme. We set the stabilization parameter~$\varepsilon=10^{-8}$. 
\par
We start the simulations of \textalpha{}-synuclein diffusion from two different initial conditions for the misfolded proteins: in the first case, they are located in the dorsal motor nucleus, namely, at the base of the brainstem; in the second one, they are located in the limbic region, just under the ventricles~\cite{adler_unified_2019}.
\par
In Figure~\ref{fig:Solution2DPark}, we report the initial conditions at time $t=0$ and the computed solutions at different times $t=5,10,15,20$ years. It can be observed that both starting with a brainstem or a limbic predominant concentration, the next stage of the propagation is associated with the presence of Lewy bodies in both regions, see Figure~\ref{fig:Solution2DPark}, $t=10$ years. Finally, the spread of the \textalpha{-}synuclein continues into the neocortical area, as visible in Figure~\ref{fig:Solution2DPark} ($t= 15,\,20$ years). This pathological progress is coherent with the medical literature findings of USSLB theory \cite{adler_unified_2019}. Medical literature describes that a limbic-predominant Lewy body pathology is more associated with early cognitive decline \cite{cersosimo_propagation_2018}. This is coherent with the earlier involvement of the neocortex resulting from our numerical simulation.
\par
\begin{figure}[t!]
    \centering
    {\includegraphics[width=0.9\textwidth]{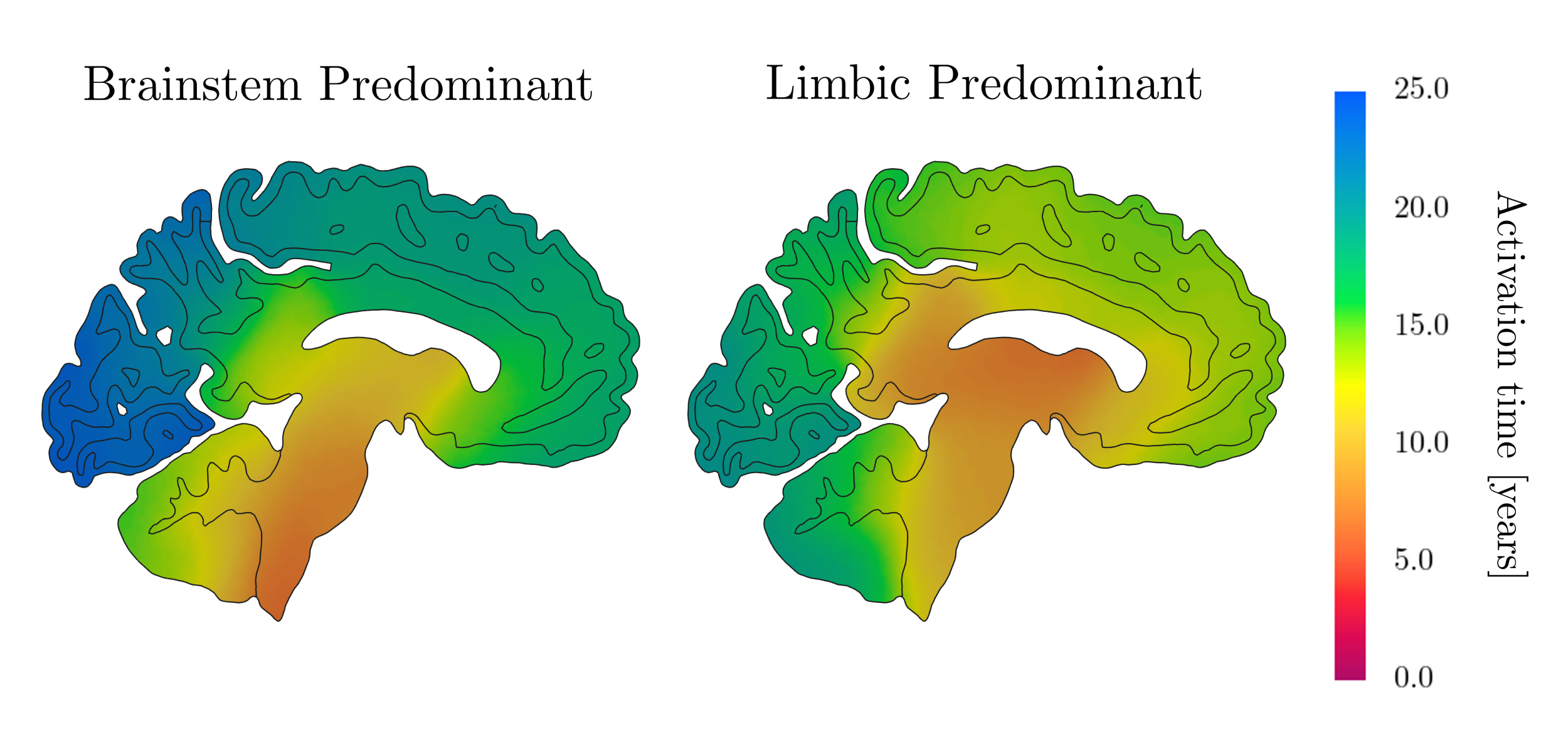}}
    \caption{Test case 3: Activation time patterns of \textalpha{}-synuclein propagation for the brainstem (left) and limbic predominant (right) prodromal evolutions of the pathology.} 
    \label{fig:ActTime2DPark}
\end{figure}
\par
Finally, to quantify the development of the pathology, we compute the activation time 
as defined in~\cite{Corti_Bonizzoni_Dede_Quarteroni_Antonietti:2023}:
\begin{equation}
\label{eq:acttime}
    \hat{t}(\boldsymbol{x}) = \int_0^T \chi_{\{c_h(\boldsymbol{x},s)< c_\mathrm{crit}\}} (\boldsymbol{x},s) \mathrm{d}s \qquad \boldsymbol{x}\in \Omega,
\end{equation}
where $\chi$ is the indicator function and $c_\mathrm{crit}$ is the critical value of the pathological protein concentration. This indicator provides a measure of the time at which the concentration of misfolded proteins exceeds a certain threshold, after which we assume the neurons in a specific region could be affected by pathological communication induced by the induced toxicity \cite{braakStagingBrainPathology2003}.
In this test case, we fix the threshold to be equal to $c_\mathrm{crit}=0.95$. We report the computed activation time in Figure \ref{fig:ActTime2DPark}. The pattern of the brainstem predominant activation is consistent with existing literature results obtained using other numerical methods~\cite{Corti_Bonizzoni_Dede_Quarteroni_Antonietti:2023,Corti_Bonizzoni_Antonietti:2024}. 
Comparing the result of the cases with brainstem and limbic predominant phase, we observe that, in the latter case, the neocortex is affected for a shorter duration (12--20 years). In contrast, the average activation time of the neocortex in the brainstem predominant case is about
18--25 years. This confirms the early cognitive decline of the limbic predominant Lewy body pathology \cite{cersosimo_propagation_2018}.

\section{Conclusions}
This work has introduced a structure-preserving, high-order, 
unconditionally stable numerical method
for approximating the solution to the Fisher-Kolmogorov equation on general polytopal meshes, with a particular focus
on modeling the spread of misfolded proteins in neurodegenerative diseases. The proposed approach is based on an entropy-variable reformulation of the model problem to ensure positivity, boundedness, and satisfaction of a discrete entropy-stability inequality at the discrete level. The numerical scheme employs a local discontinuous Galerkin (LDG) method on polytopal meshes for the spatial discretization, coupled with a $ν$-step backward differentiation formula for the time discretization.
We have discussed implementation details, highlighting the derivation and solution of the linear systems arising from Newton's iterations. The accuracy and robustness of the proposed BDF-LDG method have been demonstrated through a series of numerical experiments. We have also discussed numerical results obtained for the application of interest, namely, modeling the dynamics of neurodegenerative disorders.
In this context, we have presented numerical results showcasing the simulation of \textalpha{}-synuclein propagation in a two-dimensional brain geometry reconstructed from MRI data, thereby offering a promising computational tool for studying synucleopathies, such as Parkinson's disease and dementia with Lewy bodies.
Future developments include the theoretical analysis of the proposed approach, which is currently under investigation. Additionally, extending to three-dimensional testing represents a natural direction for further research. 
Moreover, integrating data assimilation techniques and uncertainty quantification will be crucial for enhancing the clinical relevance of the proposed computational framework.

\section*{Acknowledgements}
The first three authors were partially supported by the European Union (ERC Synergy, NEMESIS, project number
101115663). Views and opinions expressed are, however, those of the authors only and do not necessarily
reflect those of the EU or the ERC Executive Agency.
The first three authors are also members of the INdAM-GNCS group.
This research was funded in part by the Austrian Science Fund (FWF) project 10.55776/F65. The present research is part of the activities of Dipartimento di Eccellenza 2023-2027.

The brain MRI images were provided by OASIS-3: Longitudinal Multimodal Neuroimaging: Principal Investigators: T. Benzinger, D. Marcus, J. Morris; NIH P30 AG066444, P50 AG00561, P30 NS09857781, P01 AG026276, P01 AG003991, R01 AG043434, UL1 TR000448, R01 EB009352. AV-45 doses were provided by
Avid Radiopharmaceuticals, a wholly-owned subsidiary of Eli Lilly.


\end{document}